\documentclass[11pt]{article}
\usepackage{latexsym}

\font\Cp = msbm10

\newcommand{\Bbb}{\hbox{\Cp B}}
\newcommand{\Ccc}{\hbox{\Cp C}}
\newcommand{\Nnn}{\hbox{\Cp N}}
\newcommand{\Ppp}{\hbox{\Cp P}}

\newcommand{\Ttt}{\hbox{\Cp T}}

\newcommand{\hz}{\hat{0}}
\newcommand{\ho}{\hat{1}}

\newcommand{\coveredby}{\prec}

\newcommand{\qed}{\mbox{$\Box$}\vspace{\baselineskip}}

\newenvironment{proof}{\noindent {\bf Proof:}}{{\qed}}

\newenvironment{proof_}[1]{\noindent {\bf #1}}{{\qed}}

\newtheorem{theorem}{Theorem}[section]
\newtheorem{proposition}[theorem]{Proposition}
\newtheorem{lemma}[theorem]{Lemma}
\newtheorem{fact}[theorem]{Fact}
\newtheorem{definition}[theorem]{Definition}
\newtheorem{corollary}[theorem]{Corollary}
\newtheorem{example}[theorem]{Example}

\newcommand{\mb}{\overline{\mu}}
\newcommand{\join}{\vee}

\begin{document}

\title{Classification of the factorial functions of
       Eulerian binomial and Sheffer posets}

\author{{\sc Richard EHRENBORG
             and
             Margaret A.\ READDY
}}

\date{ $\:$ \\ \small Dedicated to Richard Stanley on the
occasion of his $60$th birthday.}

\maketitle


\begin{abstract}
We give a complete classification of  the factorial functions
of Eulerian binomial posets.  
The factorial function $B(n)$ either coincides with
$n!$, the factorial function of
the infinite Boolean algebra, or 
$2^{n-1}$, the factorial function of the infinite butterfly
poset.
We also classify the factorial functions for
Eulerian Sheffer posets.
An Eulerian Sheffer poset with 
binomial factorial function
$B(n) = n!$
has Sheffer factorial function $D(n)$
identical to that of 
the infinite Boolean algebra,
the infinite Boolean algebra with two new
coatoms inserted, or the infinite cubical poset.
Moreover, we are able to classify
the Sheffer factorial functions of
Eulerian Sheffer posets with 
binomial factorial function
$B(n) = 2^{n-1}$
as the doubling of an upside-down tree
with ranks $1$ and $2$ modified.

\vspace{2mm}
When we impose the further condition that a given Eulerian binomial
or Eulerian Sheffer poset is a lattice, this forces the poset
to be the infinite Boolean algebra
$\Bbb_X$ or the infinite cubical lattice
$\Ccc_X^{< \infty}$.
We also include several poset constructions
that have the same factorial functions as the infinite
cubical poset,
demonstrating
that classifying Eulerian Sheffer posets is a difficult
problem.
\end{abstract}

\section{Introduction}
\label{section_introduction}
\setcounter{equation}{0}

Binomial posets were introduced by Doubilet, Rota and
Stanley~\cite{D_R_S} to explain why generating functions 
naturally occurring in combinatorics have certain
forms. They are highly regular posets since the essential requirement
is that every two intervals of the same length have
the same number of maximal chains. 
As a result, many poset invariants are 
determined. For instance, the quintessential M\"obius function
is described by the generating function identity
\begin{equation}
      \sum_{n \geq 0} \mu(n) \cdot \frac{t^{n}}{B(n)}
   =
     \left(
        \sum_{n \geq 0} \frac{t^{n}}{B(n)}
     \right)^{-1}   ,  
\label{equation_generating_function}
\end{equation}
where $\mu(n)$ is the M\"obius function of an $n$-interval
and $B(n)$ is the factorial function, that is,
the number of maximal chains in an $n$-interval.
A
 binomial poset
is required to  contain an infinite chain
so that there are intervals of any length in the poset.

A graded poset is {\em Eulerian}
if its M\"obius function
is given by 
$\mu(x,y) = (-1)^{\rho(y)-\rho(x)}$
for all $x \leq y$ in the poset.
Equivalently, 
every interval of the poset satisfies the
Euler-Poincar\'e relation: the number of elements
of even rank is equal to the number of elements of odd rank
in the interval. 
The foremost example of Eulerian posets are 
face lattices of convex polytopes
and more generally, the face posets of regular $CW$-spheres.
Hence there is much geometric and topological interest
in understanding them.

A natural question arises: which binomial posets are Eulerian?
By equation~(\ref{equation_generating_function})
it is clear that 
the Eulerian property can be determined by knowing the
factorial function. In this paper we classify
the factorial functions of Eulerian binomial posets.
There are two possibilities,
namely,  for the factorial 
function to correspond to 
that of the infinite Boolean algebra
or  the infinite butterfly poset.

Notice that this classification is on the
level of the factorial function,
not the poset itself.
There are more Eulerian binomial posets than
these two essential examples.  See Examples~\ref{example_product_Q}
through~\ref{example_ideal_splitting}.
However, we are able to classify
the intervals of Eulerian binomial posets.
They are either isomorphic to the finite Boolean algebra or
the finite butterfly poset.

{\em Sheffer posets} were introduced by Reiner~\cite{Reiner} and
independently by Ehrenborg and Readdy~\cite{Ehrenborg_Readdy}.  
A Sheffer poset requires
the number of maximal chains
of an interval $[x,y]$ of length $n$ 
to be  given by
$B(n)$ if $x > \hz$ and $D(n)$ if $x = \hz$.
The upper intervals 
$[x,y]$ where $x > \hz$ have the property of being
binomial. Hence the interest is to understand the Sheffer intervals
$[\hz,y]$.
Just like binomial posets, the M\"obius function
is completely determined:
\begin{equation}
   \sum_{n \geq 1} \mb(n) \frac{t^n}{D(n)}
    =
  -
  \left(
    \sum_{n \geq 1}   \frac{t^n}{D(n)}
  \right)
\cdot 
  \left(
    \sum_{n \geq 0}   \frac{t^{n}}{B(n)}
  \right)^{-1} ,  
\end{equation}
where
$\mb$ is the M\"obius function of a Sheffer interval
of length $n$;  see~\cite{Ehrenborg_Readdy,Reiner}.

The classic example of a Sheffer poset
is the infinite cubical poset (see
Example~\ref{example_cubical}).
In this case,
every interval $[x,y]$ of length $n$, where $x$
is not the minimal element $\hz$, has $n!$ maximal chains.
In fact, every such interval is isomorphic to a Boolean algebra.
Intervals of the form $[\hz,y]$ have $2^{n-1} \cdot (n-1)!$ 
maximal chains and are isomorphic to the face
lattice of a finite dimensional cube.

In Sections~\ref{section_Eulerian_Sheffer}
and~\ref{section_Eulerian_Sheffer_with_n_factorial}
we completely classify the factorial functions
of Eulerian Sheffer posets.
The factorial function $B(n)$
follows from the classification of binomial posets.
The pair of factorial functions $B(n)$ and $D(n)$
fall into  three cases
(see Theorem~\ref{theorem_Sheffer_Boolean})
and
one infinite class
(Theorem~\ref{theorem_Sheffer_butterfly}).
Furthermore,
for the infinite class
we can describe the underlying Sheffer intervals;
see Theorem~\ref{theorem_Sheffer_butterfly_II_second_try}.
For two of the three cases in 
Theorem~\ref{theorem_Sheffer_Boolean}
we can also 
classify the Sheffer intervals.
For the third case we 
construct a multitude
of examples of Sheffer posets.
See Examples~\ref{example_not_cubical},
\ref{example_not_cubical_II},
\ref{example_not_cubical_III}
and~\ref{example_not_cubical_IV}.
It is striking that 
we can find many Sheffer posets
having the same factorial functions as the
infinite cubical lattice, but with the Sheffer intervals
not isomorphic to the finite cubical lattice.
However, once we require
each Sheffer interval to be a lattice then
we obtain that 
the Sheffer intervals are isomorphic to
cubical lattices.

When we impose the further condition that a given Eulerian binomial
or Eulerian Sheffer poset is a lattice, this forces the poset
to be the infinite Boolean algebra
$\Bbb_X$ or the infinite cubical lattice
$\Ccc_X^{< \infty}$.
See Examples~\ref{example_cardinals}
and~\ref{example_cubical_lattice_indeed}.

The classification of the factorial functions hinges
on the condition that the posets under consideration contain
an infinite chain. In the concluding remarks, we discuss
what could happen if this condition is removed.
We give examples of finite posets whose factorial functions
behave like the face lattice of the dodecahedron,
but which themselves are not isomorphic to this lattice.
This is part of a potentially large class of Eulerian posets
which are not polytopal-based.


\section{Eulerian binomial posets}
\setcounter{equation}{0}

\begin{definition}
A locally finite poset $P$ with $\hz$ is called a {\em binomial poset} 
if it satisfies the
following three conditions:
\begin{enumerate}
\item[(i)] $P$ contains an infinite chain.
\item[(ii)] Every interval $[x,y]$ is graded;
      hence $P$ has rank function $\rho$.
      If $\rho(x,y) = n$, then we call $[x,y]$ an {\em $n$-interval.}
\item[(iii)] For all $n \in \Nnn$, any two $n$-intervals contain the same
      number $B(n)$ of maximal chains. We call $B(n)$ the
      {\em factorial function} of $P$.
\end{enumerate}
If $P$ does not satisfy condition~(i)
and has a unique maximal element
then we say $P$ is a {\em finite binomial poset}.
\label{definition_binomial}
\end{definition}

For standard poset terminology, we refer the reader to \cite{Stanley_b}.
The number of elements of rank $k$ in an $n$-interval is given
by $B(n)/(B(k) \cdot B(n-k))$. In particular,
an $n$-interval has $A(n) = B(n)/B(n-1)$ atoms (and coatoms).
The function $A(n)$ is called the {\em atom function}
and expresses the factorial function as
$B(n) = A(n) \cdot A(n-1) \cdots A(1)$.
Directly we have $B(0) = B(1) = A(1) = 1$.
Since the atoms of an $(n-1)$-interval are contained
among the set of atoms of an $n$-interval,
the inequality $A(n-1) \leq A(n)$ holds.
Observe that if a finite binomial poset has rank $j$, 
the factorial and atom functions are
only defined up to $j$.
For further background 
material on binomial posets,
see~\cite{D_R_S,Stanley_a,Stanley_b}.

\begin{example}
{\rm
Let $\Bbb$ be the collection of finite subsets
of the positive integers ordered by inclusion.
The poset $\Bbb$ is a binomial poset with
factorial function $B(n) = n!$ and atom function
$A(n) = n$. An $n$-interval is isomorphic
to the Boolean algebra $B_{n}$.
This example is the {\em infinite Boolean algebra}.
}
\end{example}

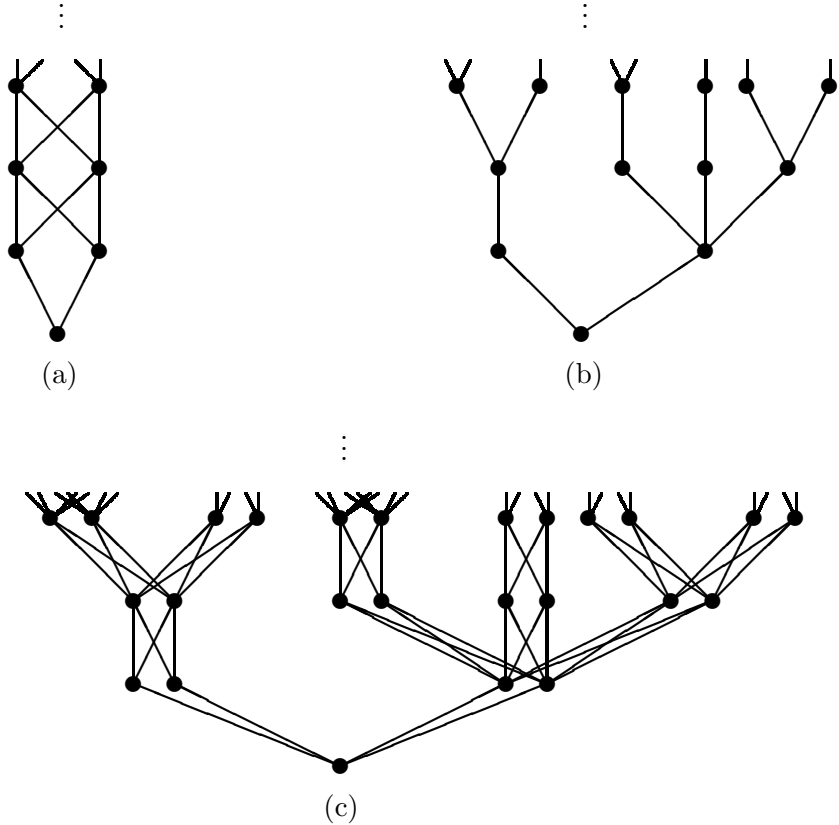
\begin{figure}
\setlength{\unitlength}{1.1mm}
\begin{center}
\begin{picture}(10,40)(0,-7)
\thicklines
\put(5,0){\line(-1,2){5}}
\put(5,0){\line(1,2){5}}
\put(5,0){\circle*{2}}

\put(0,10){\line(0,1){10}}
\put(0,10){\line(1,1){10}}
\put(0,10){\circle*{2}}

\put(10,10){\line(-1,1){10}}
\put(10,10){\line(0,1){10}}
\put(10,10){\circle*{2}}

\put(0,20){\line(0,1){10}}
\put(0,20){\line(1,1){10}}
\put(0,20){\circle*{2}}

\put(10,20){\line(-1,1){10}}
\put(10,20){\line(0,1){10}}
\put(10,20){\circle*{2}}

\qbezier(0,30)(0,32)(0,33)
\qbezier(0,30)(2,32)(3,33)
\put(0,30){\circle*{2}}

\qbezier(10,30)(8,32)(7,33)
\qbezier(10,30)(10,32)(10,33)
\put(10,30){\circle*{2}}

\put(5,37){$\vdots$}
\put(3,-6){(a)}
\end{picture}
\hspace*{45 mm}
\begin{picture}(45,40)(35,-7)
\thicklines
\put(50,0){\line(-1,1){10}}
\put(50,0){\line(3,2){15}}
\put(50,0){\circle*{2}}

\put(40,10){\line(0,1){10}}
\put(40,10){\circle*{2}}

\put(65,10){\line(-1,1){10}}
\put(65,10){\line(0,1){10}}
\put(65,10){\line(1,1){10}}
\put(65,10){\circle*{2}}

\put(40,20){\line(-1,2){5}}
\put(40,20){\line(1,2){5}}
\put(40,20){\circle*{2}}

\put(55,20){\line(0,1){10}}
\put(55,20){\circle*{2}}

\put(65,20){\line(0,1){10}}
\put(65,20){\circle*{2}}

\put(75,20){\line(-1,2){5}}
\put(75,20){\line(1,2){5}}
\put(75,20){\circle*{2}}

\qbezier(35,30)(34,32)(33.5,33)
\qbezier(35,30)(36,32)(36.5,33)
\put(35,30){\circle*{2}}

\qbezier(45,30)(45,32)(45,33)
\put(45,30){\circle*{2}}

\qbezier(55,30)(54,32)(53.5,33)
\qbezier(55,30)(56,32)(56.5,33)
\put(55,30){\circle*{2}}

\qbezier(65,30)(65,32)(65,33)
\put(65,30){\circle*{2}}

\qbezier(70,30)(70,32)(70,33)
\put(70,30){\circle*{2}}

\qbezier(80,30)(80,32)(80,33)
\put(80,30){\circle*{2}}

\put(50,37){$\vdots$}
\put(48,-6){(b)}
\end{picture}

\begin{picture}(90,45)(15,0)
\thicklines
\put(50,0){\line(-5,2){25}}
\put(50,0){\line(-2,1){20}}
\put(50,0){\line(2,1){20}}
\put(50,0){\line(5,2){25}}
\put(50,0){\circle*{2}}

\put(25,10){\line(0,1){10}}
\put(25,10){\line(1,2){5}}
\put(25,10){\circle*{2}}

\put(30,10){\line(-1,2){5}}
\put(30,10){\line(0,1){10}}
\put(30,10){\circle*{2}}

\put(70,10){\line(-2,1){20}}
\put(70,10){\line(-3,2){15}}
\put(70,10){\line(0,1){10}}
\put(70,10){\line(1,2){5}}
\put(70,10){\line(2,1){20}}
\put(70,10){\line(5,2){25}}
\put(70,10){\circle*{2}}

\put(75,10){\line(-5,2){25}}
\put(75,10){\line(-2,1){20}}
\put(75,10){\line(-1,2){5}}
\put(75,10){\line(0,1){10}}
\put(75,10){\line(3,2){15}}
\put(75,10){\line(2,1){20}}
\put(75,10){\circle*{2}}

\put(25,20){\line(-1,1){10}}
\put(25,20){\line(-1,2){5}}
\put(25,20){\line(1,1){10}}
\put(25,20){\line(3,2){15}}
\put(25,20){\circle*{2}}

\put(30,20){\line(-3,2){15}}
\put(30,20){\line(-1,1){10}}
\put(30,20){\line(1,2){5}}
\put(30,20){\line(1,1){10}}
\put(30,20){\circle*{2}}

\put(50,20){\line(0,1){10}}
\put(50,20){\line(1,2){5}}
\put(50,20){\circle*{2}}

\put(55,20){\line(-1,2){5}}
\put(55,20){\line(0,1){10}}
\put(55,20){\circle*{2}}

\put(70,20){\line(0,1){10}}
\put(70,20){\line(1,2){5}}
\put(70,20){\circle*{2}}

\put(75,20){\line(-1,2){5}}
\put(75,20){\line(0,1){10}}
\put(75,20){\circle*{2}}

\put(90,20){\line(-1,1){10}}
\put(90,20){\line(-1,2){5}}
\put(90,20){\line(1,1){10}}
\put(90,20){\line(3,2){15}}
\put(90,20){\circle*{2}}

\put(95,20){\line(-3,2){15}}
\put(95,20){\line(-1,1){10}}
\put(95,20){\line(1,2){5}}
\put(95,20){\line(1,1){10}}
\put(95,20){\circle*{2}}

\qbezier(15,30)(13,32)(12,33)
\qbezier(15,30)(14,32)(13.5,33)
\qbezier(15,30)(17,32)(18,33)
\qbezier(15,30)(18,32)(19.5,33)
\put(15,30){\circle*{2}}

\qbezier(20,30)(17,32)(15.5,33)
\qbezier(20,30)(18,32)(17,33)
\qbezier(20,30)(21,32)(21.5,33)
\qbezier(20,30)(22,32)(23,33)
\put(20,30){\circle*{2}}

\qbezier(35,30)(35,32)(35,33)
\qbezier(35,30)(36,32)(36.5,33)
\put(35,30){\circle*{2}}

\qbezier(40,30)(39,32)(38.5,33)
\qbezier(40,30)(40,32)(40,33)
\put(40,30){\circle*{2}}

\qbezier(50,30)(48,32)(47,33)
\qbezier(50,30)(49,32)(48.5,33)
\qbezier(50,30)(52,32)(53.5,33)
\qbezier(50,30)(53,32)(54.5,33)
\put(50,30){\circle*{2}}

\qbezier(55,30)(52,32)(50.5,33)
\qbezier(55,30)(53,32)(52,33)
\qbezier(55,30)(56,32)(56.5,33)
\qbezier(55,30)(57,32)(58,33)
\put(55,30){\circle*{2}}

\qbezier(70,30)(70,32)(70,33)
\qbezier(70,30)(71,32)(71.5,33)
\put(70,30){\circle*{2}}

\qbezier(75,30)(74,32)(73.5,33)
\qbezier(75,30)(75,32)(75,33)
\put(75,30){\circle*{2}}

\qbezier(80,30)(80,32)(80,33)
\qbezier(80,30)(81,32)(81.5,33)
\put(80,30){\circle*{2}}

\qbezier(85,30)(84,32)(83.5,33)
\qbezier(85,30)(85,32)(85,33)
\put(85,30){\circle*{2}}

\qbezier(100,30)(100,32)(100,33)
\qbezier(100,30)(101,32)(101.5,33)
\put(100,30){\circle*{2}}

\qbezier(105,30)(104,32)(103.5,33)
\qbezier(105,30)(105,32)(105,33)
\put(105,30){\circle*{2}}

\put(50,37){$\vdots$}
\put(48,-6){(c)}
\end{picture}
\end{center}
\caption{(a) the infinite butterfly poset $\Ttt$,
         (b) an infinite tree $Q$,
         (c) and the rank product $\Ttt * Q$,
which has the same factorial function as
the butterfly poset.}
\label{figure_butterfly_tree}
\end{figure}

\begin{example}
{\rm
Let $\Ttt$ be the {\em infinite butterfly poset}, that is,
$\Ttt$ consists of the elements  
$\{\hz\} \cup \left(\Ppp \times \{1,2\}\right)$ 
where 
$(n,i) \prec (n+1,j)$ for all $i,j \in \{1,2\}$
and $\hz$ is the unique minimal element;
see Figure~\ref{figure_butterfly_tree}~(a).
The poset $\Ttt$ is a binomial poset.  It has
factorial function $B(n) = 2^{n-1}$ for $n \geq 1$
and atom function
$A(n) = 2$ for $n \geq 2$.
Let $T_{n}$ denote an $n$-interval in $\Ttt$.
}
\end{example}

\begin{example}
{\rm
Given two ranked posets $P$ and $Q$, define the {\em rank product} $P*Q$
by
$$ P * Q
  =
   \{(x,z) \in P \times Q \:\: : \:\: \rho_{P}(x) = \rho_{Q}(z) \} .  $$
Define the order relation by
$(x,y) \leq_{P * Q} (z,w)$ if
$x \leq_{P} z$ and $y \leq_{Q} w$.
If $P$ and $Q$ are binomial posets then so is
the poset $P*Q$.  It has the factorial function
$B_{P*Q}(n) = B_{P}(n) \cdot B_{Q}(n)$.
This example is due to Stanley~\cite[Example~3.15.3~d]{Stanley_b}.
The rank product is also known as the Segre product;
see~\cite{Bjorner_Welker}
}
\end{example}

\begin{example}
{\rm
For $q \geq 2$ let $P_{q}$ be the face poset of an $q$-gon.
Observe that this is a finite binomial poset of rank $3$
with the factorial function $B(2) = 2$ and $B(3) = 2q$.
Let $q_{1}, \ldots, q_{r}$ be a list of integers 
with each $q_{i} \geq 2$.
Let $P_{q_{1}, \ldots, q_{r}}$ be the poset obtained by
identifying all the minimal elements of $P_{q_{1}}$ through
$P_{q_{r}}$ and identifying all the maximal elements.
This is also a binomial poset with factorial function
$B(2) = 2$ and $B(3) = 2 \cdot (q_{1} + \cdots + q_{r})$.
It is straightforward to see that
each rank $3$ binomial poset with $B(2) = 2$ is of this form.
}
\label{example_rank_3}
\end{example}

A finite graded poset
is said to satisfy the {\em Euler-Poincar\'e relation}
if it 
has the same number of elements of even rank as of odd rank.
A poset is called {\em Eulerian} if every non-singleton interval
satisfies the Euler-Poincar\'e relation.
Equivalently, a poset $P$ is Eulerian if its
M\"obius function satisfies 
$\mu(x,y) = (-1)^{\rho(y) - \rho(x)}$ for all
$x \leq y$ in $P$.

\begin{lemma}
Let $P$ be a graded poset of odd rank such that every proper interval
of $P$ is Eulerian.
Then $P$ is an Eulerian poset.
\label{lemma_odd_ranked_posets}
\end{lemma}
This is Exercise 69c in~\cite{Stanley_b}.
Also this lemma is implicit in
the two papers~\cite{Billera_Liu,Ehrenborg_Readdy_homology}.
A three-line proof is as follows.

\begin{proof_}{{\bf Proof of Lemma~\ref{lemma_odd_ranked_posets}:}}
We know the M\"obius function of $P$ satisfies
$\mu(x,y) = (-1)^{\rho(y) - \rho(x)}$
for $\rho(y) - \rho(x) \leq n-1$,
where $n$ is the rank of $P$.
Now
$1 + \mu(\hz,\ho)  =  - \sum_{\hz < x < \ho} (-1)^{\rho(x)}
                   =    \sum_{\hz < x < \ho} (-1)^{n - \rho(x)}
                   =    -1 - \mu(\hz, \ho)$.
This yields
$\mu(\hz,\ho) = -1 = (-1)^n$, proving that $P$ is Eulerian.
\end{proof_}

We now conclude
\begin{proposition}
To verify that a poset is Eulerian it is enough to verify
that every interval of even rank satisfies the
Euler-Poincar\'e relation.
\label{proposition_only_even_ranks}
\end{proposition}

For an $n$-interval of an Eulerian binomial poset
the Euler-Poincar\'e relation states 
\begin{equation}
      \sum_{k=0}^{n}
          (-1)^{k}
             \cdot
          \frac{B(n)}{B(k) \cdot B(n-k)}   
   =
      0        .
\label{equation_Euler_Poincare_binomial}
\end{equation}
When $n$ is even, it follows from
(\ref{equation_Euler_Poincare_binomial}) that
$B(n)$ is determined by 
$B(0), B(1), \ldots, B(n-1)$.
Also observe that $B(2) = A(2) = 2$ since every $2$-interval
is a diamond.

\begin{theorem}
Let $P$ be an Eulerian binomial poset
with factorial function $B(n)$.
Then either 
\begin{enumerate}
\item[(i)]
the factorial function $B(n)$ is given by $B(n) = n!$
and every $n$-interval is isomorphic to the Boolean algebra $B_{n}$,
or 

\item[(ii)]
the factorial function $B(n)$ is given by 
$B(0) = 1$ and $B(n) = 2^{n-1}$
and every $n$-interval is isomorphic
to the butterfly poset $T_{n}$.
\end{enumerate}
\label{theorem_binomial}
\end{theorem}

It is tempting to state this theorem
as, ``There are only two Eulerian binomial posets,
namely, the infinite Boolean algebra $\Bbb$
and the infinite butterfly poset $\Ttt$.''
However, this is false.  The next three examples
demonstrate this.
\begin{example}
{\rm
Let $Q$
be an infinite poset with a minimal element $\hz$
containing an infinite chain
such that every interval of the form $[\hz, x]$ is a chain.
Observe the poset $Q$ is an infinite tree and, in fact, 
is a binomial poset with factorial function $B(n) = 1$.
Thus we know that both $\Bbb * Q$ and $\Ttt * Q$
are Eulerian binomial posets.
See Figure~\ref{figure_butterfly_tree} for an example.
When the poset $Q$ is different from an infinite
chain, we have that
$\Bbb * Q \not\cong \Bbb$ and
$\Ttt * Q \not\cong \Ttt$.
This follows since
in the two posets 
$\Bbb$ and $\Ttt$ every pair of elements
has an upper bound, that is, the two posets are
confluent.
This property does not hold in the tree $Q$ and hence
not in the rank products
$\Bbb * Q$ and $\Ttt * Q$ 
either.
}
\label{example_product_Q}
\end{example}

\begin{example}
{\rm
For each infinite cardinal $\kappa$ there is a Boolean algebra
consisting of all finite subsets of a set $X$ of cardinality $\kappa$.
We denote this poset by $\Bbb_X$.
Observe that different cardinals give rise to 
non-isomorphic Boolean algebras.
}
\label{example_cardinals}
\end{example}

\begin{example}
{\rm
Let $P$ be a binomial poset and $I$ a nonempty lower order ideal of
$P$. Construct a new poset by taking the Cartesian product
of the poset $P$ with the two element antichain $\{a,b\}$,
and identify elements of the form $(x,a)$ and $(x,b)$ if
$x$ lies in the ideal $I$. The new poset is also binomial
and has the same factorial function as $P$.
}
\label{example_ideal_splitting}
\end{example}


We now state a very useful lemma.

\begin{lemma}
Let $P$ and $P^{\prime}$ be two Eulerian binomial posets
having atom functions $A(n)$ and $A^{\prime}(n)$ 
which agree
for $n \leq 2m$, where $m \geq 2$.
Then the following equality holds:
\begin{equation}
        \frac{1}{A(2m+1)}
     \cdot
        \left(
          1 - \frac{1}{A(2m+2)}
        \right) 
    =
        \frac{1}{A^{\prime}(2m+1)}
     \cdot
        \left(
          1 - \frac{1}{A^{\prime}(2m+2)}
        \right) 
   .  
\label{equation_A_A}
\end{equation}
\label{lemma_A_A}
\end{lemma}
\begin{proof}
Let $B(n)$ and $B^{\prime}(n)$ 
be the factorial functions for $P$, respectively $P^{\prime}$.
By the Euler-Poincar\'e relation, we have
$$
\sum_{k=0}^{2m+2}
   (-1)^{k} \cdot
   \frac{1}{B(k) \cdot B(2m+2-k)} 
          = 0 =
\sum_{k=0}^{2m+2}
   (-1)^{k} \cdot
   \frac{1}{B^{\prime}(k) \cdot B^{\prime}(2m+2-k)} 
.
$$
Cancelling all the terms
where $B$ and $B^{\prime}$
agree, we have
\begin{eqnarray*}
&&
      \frac{2}{A(2m+2) \cdot A(2m+1) \cdot B(2m)} 
     -
        \frac{2}{A(2m+1) \cdot B(2m)}  \\
&=&
        \frac{2}{A^{\prime}(2m+2) \cdot A^{\prime}(2m+1) \cdot B(2m)} 
     -
        \frac{2}{A^{\prime}(2m+1) \cdot B(2m)}     . 
\end{eqnarray*}
Cancelling common factors, we obtain the desired equality.
\end{proof}

As a corollary to Lemma~\ref{lemma_A_A} we have:

\begin{corollary}
Let $P$ and $P^{\prime}$
be two Eulerian binomial posets satisfying the conditions
in Lemma~\ref{lemma_A_A}.
Assume furthermore there is  a lower and an upper bound
for $A^{\prime}(2m+2)$ of the form
$L \leq A^{\prime}(2m+2) < U$. 
Let $x$ be the left-hand side of equation~(\ref{equation_A_A}).
Then we obtain 
a lower and an upper bound for
$A^{\prime}(2m+1)$, namely
\begin{equation}
     \frac{1}{x} \cdot \left(1 - \frac{1}{L}\right)
   \leq
     A^{\prime}(2m+1)
   <
     \frac{1}{x} \cdot \left(1 - \frac{1}{U}\right) .  
\label{equation_A_A_bound}
\end{equation}
\label{corollary_A_A_bound}
\end{corollary}
We see that these bounds can be improved by using that
$A^{\prime}(2m+1)$ is in fact an integer.

\begin{proposition}
Let $P^{\prime}$ be an Eulerian binomial poset
with factorial function $B^{\prime}(n)$
satisfying
$B^{\prime}(3) = 6$.
Then the factorial function
is given by $B^{\prime}(n) = n!$.
\label{proposition_Boolean}
\end{proposition}
\begin{proof}
Let $P$ be the infinite Boolean algebra $\Bbb$ with atom function
$A(n) = n$ and factorial function $B(n) = n!$.
We will prove that the two factorial functions
$B(n)$ and $B^{\prime}(n)$ are identical, equivalently
that the two atom functions
$A(n)$ and $A^{\prime}(n)$ are equal.

Assume that the two atom functions $A$ and $A^{\prime}$
agree up to $2m = j$.
Since $A(n) = n$ the left-hand side of
equation~(\ref{equation_A_A})
is equal to $1/(j+2)$.
We have the following bounds for $A^{\prime}(j+2)$:
$j = A^{\prime}(j) \leq A^{\prime}(j+2) < \infty$.
Applying Corollary~\ref{corollary_A_A_bound} we obtain the following
bounds on $A^{\prime}(j+1)$:
$$     j+1-\frac{2}{j} \leq A^{\prime}(j+1) < j+2  .  $$
Since $A^{\prime}(j+1)$ is an integer and $j \geq 4$
we conclude that
$A^{\prime}(j+1) = j+1$.
This implies that
$A^{\prime}(j+2) = j+2$
and hence 
we conclude the two atom functions are equal.
\end{proof}

\begin{proposition}
Let $P$ be a finite binomial poset
of rank $n$ with factorial function $B(k) = k!$
for $k \leq n$.
Then the poset $P$
is isomorphic to the Boolean algebra $B_{n}$.
\label{proposition_Boolean_structure}
\end{proposition}
\begin{proof}
Directly the result  is true for $n \leq 2$.
Assume it is true for all posets of rank $n-1$ and consider a
poset $P$ of rank $n$.
Since $P$ is a binomial poset 
with factorial function $B(k) = k!$,
we know that
the number of elements of rank $k$ in $P$
is given by ${n \choose k}$. Especially, the cardinality
of $P$ is given by $2^{n}$.
Let~$c$ be a coatom in the poset.
Observe that the interval $[\hz,c]$ is isomorphic to $B_{n-1}$
by the induction hypothesis and hence the coatom $c$ is greater
than all but one atom $a$ in the poset $P$. Similarly,
the interval $[a,\ho]$ is also isomorphic to $B_{n-1}$.
Since the two intervals $[a,\ho]$ and $[\hz,c]$ are disjoint
and have the same cardinality $2^{n-1}$, the poset
$P$ is the disjoint union of these two intervals.

Using the binomial property of $P$, an element
$z$ of rank $k$ in the lower interval $[\hz,c]$ is covered
by $n-k$ elements in
the poset $P$ and by $n-k-1$ elements in
the interval $[\hz,c]$. 
Thus there is a unique element in $[a,\ho]$ that covers $z$.
Denote this element by $\varphi(z)$. 
By a similar argument we obtain that
$\varphi$ is a bijective function from $[\hz,c]$ to $[a,\ho]$.
Let $z \prec w$ be a cover relation in $[\hz,c]$.
Consider the $2$-interval $[z,\varphi(w)]$. 
As every $2$-interval is
a diamond there is an element $v$ different from $w$ such that
$z \prec v \prec \varphi(w)$. Since $w$ is the unique element in $[\hz,c]$
that is covered by $\varphi(w)$,
the element $v$ belongs to the upper interval
$[a,\ho]$. Also,
the element  $\varphi(z)$ is the unique element in the upper interval
that covers $z$, we conclude that $v = \varphi(z)$ and especially
$\varphi(w)$ covers $\varphi(z)$.
Hence the function $\varphi$ is order-preserving. By the symmetric argument
$\varphi^{-1}$ is also order-preserving. 
Therefore the poset $P$ is
the Cartesian product of $[\hz,c]$ with the two element poset $B_{1}$
and we conclude that $P$ is isomorphic to the Boolean algebra
$B_{n}$.
\end{proof}

\begin{proposition}
Let $P^{\prime}$ be an Eulerian binomial poset
with factorial function $B^{\prime}(n)$
satisfying
$B^{\prime}(3) = 4$.
Then the factorial function
is given by $B^{\prime}(n) = 2^{n-1}$ for $n \geq 1$.
\label{proposition_butterfly}
\end{proposition}
\begin{proof}
Let $P$ be the butterfly poset $\Ttt$
and $A(n)$ its atom function, where
$A(1) = 1$ and $A(n) = 2$ for $n \geq 2$.
Similar to the proof of
Proposition~\ref{proposition_Boolean}
we consider how $A(n)$ and $A^{\prime}(n)$ relate.

Assume that the two atom functions agree up to $2m = j$.
Now the left-hand side of equation~(\ref{equation_A_A})
is equal to $1/4$.
For $A^{\prime}(j+2)$ we have the bounds
$2 = A^{\prime}(j) \leq A^{\prime}(j+2) < \infty$.
Applying Corollary~\ref{corollary_A_A_bound} we obtain 
$$     2 \leq A^{\prime}(j+1) < 4   .  $$

Consider now the possibility that
$A^{\prime}(j+1) = 3$.
Let $[x,y]$ be a $(j+1)$-interval in $P^{\prime}$.
For $1 \leq k \leq j$ there are
$B^{\prime}(j+1)/(B^{\prime}(k) \cdot B^{\prime}(j+1-k))
  =
 3 \cdot 2^{j-1} / (2^{k-1} \cdot 2^{j-k}) = 3$
elements of rank $k$ in this interval.
Let~$c$ be a coatom. The interval $[x,c]$ has
two atoms, say $a_{1}$ and $a_{2}$.
Moreover, the interval $[x,c]$ has
two elements of rank $2$, say $b_{1}$ and $b_{2}$.
Moreover we know that each $b_{j}$ covers each~$a_{i}$.
Let~$a_{3}$ and~$b_{3}$ be the third atom, respectively
the third rank $2$ element, in the interval $[x,y]$.
We know that~$b_{3}$ covers two atoms in $[x,y]$.
One of them must be $a_{1}$ or $a_{2}$, say $a_{1}$.
But then $a_{1}$ is covered by the three elements $b_{1}$, 
$b_{2}$ and $b_{3}$.
But this contradicts
the fact that each atom is covered by exactly
two elements.
Hence this rules out the case
$A^{\prime}(j+1) = 3$.

The only remaining possibility is
$A^{\prime}(j+1) = 2$, implying
$A^{\prime}(j+2) = 2$. Hence the atom functions
$A(n)$ and $A^{\prime}(n)$ are equal.
\end{proof}

\begin{lemma}
Let $P$ be a finite binomial poset
with factorial function $B(k) = 2^{k-1}$ for $1 \leq k \leq n$.
Then the poset $P$ is isomorphic to the butterfly poset $T_{n}$.
\label{lemma_butterfly_structure}
\end{lemma}
\begin{proof}
Directly true for $n \leq 2$.
Assume now that $n \geq 3$.
Observe that there are 
$B(n)/(B(k) \cdot B(n-k)) = 2$
elements of each rank and every element
of rank greater than or equal to $2$ covers exactly two elements.
Hence the only possibility is that the poset $P$ is isomorphic to
the butterfly poset $T_n$.
\end{proof}

\begin{proof_}{{\bf Proof of Theorem~\ref{theorem_binomial}:}}
The atom function of an Eulerian binomial poset
satisfies $2 = A(2) \leq A(3)$.
Hence $B(3) = A(3) \cdot B(2)$ is an even integer
greater than or equal to $4$.
The Euler-Poincar\'e relation implies that
$$    \frac{1}{B(4)} = \frac{1}{B(3)} - \frac{1}{8} ,  $$
implying that $B(3) < 8$. Hence there are only two 
remaining cases,
which are considered in
Propositions~\ref{proposition_Boolean}
and~\ref{proposition_butterfly}.
The corresponding structure statements are considered in
Proposition~\ref{proposition_Boolean_structure}
and Lemma~\ref{lemma_butterfly_structure}.
\end{proof_}

\begin{theorem}
Let $L$ be an Eulerian binomial poset 
which we  furthermore require to be a lattice.
Then $L$ is isomorphic to 
the Boolean algebra $\Bbb_{X}$
where $X$ is the set of atoms of the poset $L$.
\label{theorem_Eulerian_binomial_lattice}
\end{theorem}
\begin{proof}
Since every interval of $L$ is a lattice
we can rule out the butterfly factorial function.
Hence $B(n) = n!$ and every interval $[\hz,x]$
is a Boolean lattice.
Let $\varphi$ be the map from $L$ to $\Bbb_{X}$
defined by
$\varphi(x) = \{a \in X \:\: : \:\: a \leq x\}$.
The inverse of $\varphi$ is given by
${\displaystyle \varphi(Y) = \join_{a \in Y} a}$.
It is straightforward to see that
both $\varphi$ and $\varphi^{-1}$ are
order-preserving.
Hence the two lattices $L$ and $\Bbb_{X}$ are isomorphic.
\end{proof}

We end this section with a result that will be used 
in Section~\ref{section_Eulerian_Sheffer_with_n_factorial}
when we study Eulerian Sheffer posets.
\begin{proposition}
There is no finite
binomial poset $P^{\prime}$ of rank $j+1 \geq 4$ with the
atom function
$$  A^{\prime}(n)
         =
           \left\{ \begin{array}{c l} n   & \mbox{ if } n \leq j, \\
                                      j+2 & \mbox{ if } n = j+1.
                   \end{array} \right. $$
\label{proposition_special_case}
\end{proposition}
\begin{proof}
Assume that the poset $P^{\prime}$ exists. Then it has $j+2$ atoms
and
$j+2$ coatoms. Each atom $x$ lies below exactly $j$ coatoms and
each coatom $c$ lies  above exactly $j$ atoms.
Moreover, by the proof of
Proposition~\ref{proposition_Boolean}
we know that each of the intervals $[\hz,c]$ and $[x,\ho]$
is isomorphic to $B_{j}$.

Define a multigraph $G$ with the $j+2$ atoms as the vertices.
For each coatom $c$ let there be an edge~$xy$ between the two unique
atoms $x$ and $y$ such that $x,y \not\leq c$.
Since each atom is not below exactly two coatoms, each vertex of
the graph has degree equal to $2$. Hence the graph is a disjoint
union of cycles.

Pick a coatom $c$ that corresponds to an edge $xy$.
The coatom $c$ is greater than
the $j$ atoms $z_{1}, \ldots, z_{j}$.
Using that the interval $[\hz,c]$ is a Boolean algebra,
let $w_{i}$ be the unique coatom in the interval
$[\hz,c]$ that is not greater than $z_{i}$.
Let $d_{i}$ be the atom in the interval $[w_{i},\ho] \cong B_{2}$
distinct from~$c$.
Observe for $i \neq k$ we have $z_{i} < w_{k} < d_{k}$.
Hence the $j$ coatoms
$c, d_{1}, \ldots, \widehat{d_{i}}, \ldots, d_{j}$ are
all the coatoms greater than $z_{i}$.
Moreover, since $j \geq 3$ we conclude that
$d_{1}, \ldots, d_{j}$ are all distinct.

Consider the $j$ atoms below $d_{k}$. They are
$z_{1}, \ldots, \widehat{z_{k}}, \ldots, z_{j}$
and exactly one of $x$ and $y$.
Thus the edge $e_k$ corresponding to $d_{k}$ intersects
the edge $xy$. This holds for all $j$ edges $e_k$.
Hence we obtain the contradiction
$4 = \deg(x) + \deg(y) \geq 2 + j$.
Thus there is no such finite binomial poset.
\end{proof}

\section{Eulerian Sheffer posets}
\label{section_Eulerian_Sheffer}
\setcounter{equation}{0}

{\em Sheffer posets}, also know as upper binomial posets,
were first defined by Reiner~\cite{Reiner}
and independently discovered by 
Ehrenborg and Readdy~\cite{Ehrenborg_Readdy}.

\begin{definition}
A locally finite poset $P$ with $\hz$
is called a {\em Sheffer poset} if it satisfies the
following four conditions:
\begin{enumerate}
\item[(i)] $P$ contains an infinite chain.
\item[(ii)] Every interval $[x,y]$ is graded;
      hence $P$ has rank function $\rho$.
      If $\rho(x,y) = n$, then we call $[x,y]$ an {\em $n$-interval.}
\item[(iii)] Two $n$-intervals $[\hz,y]$ and $[\hz,v]$, such that
      $y \neq \hz$, $v \neq \hz$,
      have the same number $D(n)$ of maximal
      chains.
\item[(iv)] Two $n$-intervals $[x,y]$ and $[u,v]$, such that
      $x \neq \hz$, $u \neq \hz$, 
      have the same number $B(n)$ of maximal
      chains.
\end{enumerate}
As in the finite
binomial poset case,
if $P$ does not satisfy condition~(i)
and has a unique maximal element
then we say $P$ is a {\em finite Sheffer poset}.
\end{definition}
An interval of the form $[\hz,y]$ is called
a {\em Sheffer interval}, whereas an interval $[x,y]$, where
$x > \hz$, is called a {\em binomial interval}.
Similarly, the functions
$B(n)$ and $D(n)$ are called
the {\em binomial} and  {\em Sheffer factorial functions}
respectively.
The number of elements of rank $k \geq 1$
in a Sheffer interval of length~$n$ is given
by $D(n)/(D(k) \cdot B(n-k))$. Especially,
a Sheffer interval $[\hz,y]$ has $C(n) = D(n)/D(n-1)$ coatoms.
The function $C(n)$ is called the {\em coatom function}
and we have
$D(n) = C(n) \cdot C(n-1) \cdots C(1)$.
Observe that $D(1) = C(1) = 1$.

\begin{example}
{\rm
Every binomial poset is a Sheffer poset. The factorial
functions are equal, that is, $D(n) = B(n)$ for $n \geq 1$.
}
\end{example}

\begin{example}
{\rm
The rank product $P*Q$ of two Sheffer posets $P$ and $Q$
is also a Sheffer poset with
the factorial functions
$B_{P*Q}(n) = B_{P}(n) \cdot B_{Q}(n)$
and
$D_{P*Q}(n) = D_{P}(n) \cdot D_{Q}(n)$.
}
\end{example}

\begin{example}
{\rm
For a poset $P$ with a unique minimal element $\hz$,
let the {\em dual suspension}
$\Sigma^{*}(P)$ be the poset $P$ with two new elements
$a_{1}$ and $a_{2}$.
Let the order relations be as follows:
$\hz <_{\Sigma^{*}(P)} a_{i} <_{\Sigma^{*}(P)} y$
for all $y > \hz$ in $P$ and $i=1,2$.
That is, the elements $a_{1}$ and $a_{2}$
are inserted between $\hz$ and the atoms of $P$.
Clearly if $P$ is Eulerian
then so is $\Sigma^{*}(P)$. Moreover, if $P$ is a binomial poset
then $\Sigma^{*}(P)$ is a Sheffer poset with the
factorial function $D_{\Sigma^{*}(P)}(n) = 2 \cdot B(n-1)$
for $n \geq 2$.
}
\end{example}

One may extend the dual suspension $\Sigma^{*}$ by inserting
$k$ new atoms instead of $2$. Yet again it will take a binomial poset
to a Sheffer poset. However we have no need
of this extension since it does not preserve the Eulerian property.

\begin{example}
{\rm
Let $P$ be the three element poset
\setlength{\unitlength}{0.5 mm}
\begin{picture}(13,10)(-1,0)
\put(0,0){\circle*{2}}
\put(10,0){\circle*{2}}
\put(5,10){\circle*{2}}
\put(0,0){\line(1,2){5}}
\put(10,0){\line(-1,2){5}}
\put(-4,0){${\scriptstyle 0}$}
\put(11,0){${\scriptstyle 1}$}
\put(6,10){${\scriptstyle *}$}
\end{picture}.
The poset
$C_n = P^n \cup \{\hz\}$
is the 
face lattice of the $n$-dimensional cube,
also  known as the 
{\em cubical lattice}.
It is a  finite Sheffer poset
with factorial functions
$B(k) = k!$ for $k \leq n$ and $D(k) = 2^{k-1} \cdot (k-1)!$
for $1 \leq k \leq n+1$.
}
\label{example_finite_cubical}
\end{example}

For a ranked poset $P$ (not necessarily having a unique minimal element)
and a possibly infinite set $X$
define the {\em power poset} $P^{X}$ as follows.
Let the underlying set be given by
$$   P^{X}
  = 
     \left\{ f : X \rightarrow P
                     \:\: : \:\:
            \sum_{x \in X} \rho(f(x)) < \infty \right\} $$
and define the order relation by componentwise comparison,
that is, $f \leq_{P^{X}} g$ if $f(x) \leq g(x)$ for all $x$ in $X$.

\begin{example}
{\rm
Let $P$ be as in the previous example
and let $X$ be an infinite set.
The poset $\Ccc_{X} = P^{X} \cup \{\hz\}$,
that is, the poset $P^{X}$ with a new minimal element adjoined,
is a Sheffer poset.
This example is precisely the {\em infinite cubical poset} with the factorial
functions $B(n) = n!$ and $D(n) = 2^{n-1} \cdot (n-1)!$.
Similar to
Example~\ref{example_cardinals},
for different infinite cardinalities of $X$
we obtain non-isomorphic cubical posets.
Note, however, this poset is not a lattice since
the two atoms $(0,0,\ldots)$ and $(1,1,\ldots)$
do not have a join.
A Sheffer $n$-interval is isomorphic to the cubical
lattice $C_{n-1}$.
Hence, every interval in the poset 
$\Ccc_{X}$
is Eulerian.
}
\label{example_cubical}
\end{example}

\begin{example}
{\rm
Let $E_{2}, E_{3}, \ldots$ be an infinite sequence of disjoint
nonempty finite sets,
where $E_{n}$ has cardinality $e_{n}$.
Consider the poset
$$   U_{e_{2},e_{3}, \ldots}
  =
     \{\hz\} \:\: \cup \:\:
   \bigcup_{n \geq 2} \: \prod_{i \geq n} E_i  , $$
where $\prod$ stands for Cartesian product.
We make this into a ranked poset by letting $\hz$ be the minimal
element, and defining the cover relation by
$$           (x_{n}, x_{n+1}, x_{n+2}, \ldots)
      \prec         (x_{n+1}, x_{n+2}, \ldots)    ,  $$
where $x_{i} \in E_i$.
Thus the elements of $\prod_{i \geq n} E_i$ have rank $n-1$.
This poset is a Sheffer poset with the atom function
$A(n) = 1$ and coatom function is given by
$C(n) = e_{n}$ for $n \geq 2$.
We may view this poset as an ``upside-down tree''
with a minimal element attached.
}
\end{example}

Naturally, the previous example is not an Eulerian poset.
However, we can use it to construct Eulerian Sheffer posets
as the next two examples illustrate.

\begin{example}
{\rm
Recalling that $\Ttt$ denotes the infinite butterfly
poset,
consider the poset
$\Ttt * U_{e_{2},e_{3},\ldots}$,
where $e_{2} = e_{4} = e_{6} = \cdots = 1$.
This poset has the factorial
functions $B(n) = 2^{n-1}$
and
$D(n) = 2^{n-1} \cdot \prod_{i=2}^{n} e_{i}$.
In Theorem~\ref{theorem_Sheffer_butterfly}
we will observe that
the condition that $e_{2j} = 1$
implies that the poset is Eulerian.
}
\label{example_Sheffer_butterfly}
\end{example}
In general the rank product $\Ttt * P$ can be viewed as the
``doubling''
of the poset $P$. This notion was introduced by
Bayer and Hetyei in~\cite{Bayer_Hetyei}.

\begin{example}
{\rm
Let $\Bbb \cup \{\hz\}$ be the infinite Boolean algebra
with a new minimal element adjoined.
This is a Sheffer poset with factorial functions
$B(n) = n!$ and $D(n) = (n-1)!$.
Now consider the rank product
$(\Bbb \cup \{\hz\}) * U_{2,2,\ldots}$.
It has the factorial
functions $B(n) = n!$
and
$D(n) = 2^{n-1} \cdot (n-1)!$.
This poset has the same
factorial functions as the infinite cubical poset
and hence it is an Eulerian poset.
}
\label{example_not_cubical}
\end{example}

For an Eulerian Sheffer poset of rank $n$,
the Euler-Poincar\'e relation states
\begin{equation}
     1
         +
      \sum_{k=1}^{n}
          (-1)^{k}
             \cdot
          \frac{D(n)}{D(k) \cdot B(n-k)}   
   =
      0        .  
\label{equation_Euler_Poincare_Sheffer}
\end{equation}
Again by
Proposition~\ref{proposition_only_even_ranks}
this relation will only give us information for $n$ even. 
When $n = 2m$ we can write this relation as:
\begin{equation}
     \frac{2}{D(2m)}
         +
      \sum_{k=1}^{2m-1}
          (-1)^{k}
             \cdot
          \frac{1}{D(k) \cdot B(2m-k)}   
   =
      0        .  
\label{equation_Euler_Poincare_Sheffer_even}
\end{equation}
Also note that $D(2) = C(2) = 2$.

We will be using the following two facts
to exclude possible factorial functions.
\begin{fact}
{\rm 
\begin{itemize}
\item[(a)]
The inequality $A(n-1) \leq C(n) < \infty$
holds
since the set of coatoms in a Sheffer interval of rank $n$,
say $[\hz,y]$,
contains the set of coatoms in an $(n-1)$-interval $[x,y]$,
and there are a finite number of them.

\item[(b)]
The value $B(k)$ divides $C(n) \cdots C(n-k+1)$
for $n > k$,
since the number of elements of rank $n-k$
in a Sheffer interval of length $n$
is given by
the integer
$D(n)/(D(n-k) \cdot B(k)) = C(n) \cdots C(n-k+1)/B(k)$.
\end{itemize}
}
\label{fact_a_b}
\end{fact}

We end this section by classifying all Eulerian Sheffer posets
with binomial factorial function $B(n) = 2^{n-1}$. 
Theorem~\ref{theorem_Sheffer_butterfly}
classifies the Sheffer factorial function $D(n)$,
equivalently the coatom function $C(n)$, whereas 
Theorem~\ref{theorem_Sheffer_butterfly_II_second_try}
describes the Sheffer intervals.
It is noteworthy 
that Sheffer intervals in these posets are almost determined
by the factorial function $D(n)$. The Sheffer interval of rank $3$
are rather flexible
within the Sheffer and Eulerian conditions.
See Example~\ref{example_rank_3}. However, 
for higher ranks 
the structure
is then determined by the factorial function.

\begin{theorem}
Let $P$ be an Eulerian Sheffer poset
with the binomial factorial function
satisfying
$B(0) = 1$ and
$B(n) = 2^{n-1}$ for $n \geq 1$.
Then the coatom function $C(n)$ and the poset $P$
satisfy:
\begin{itemize}
\item[(i)] $C(3) \geq 2$,
and a length $3$ Sheffer interval
is isomorphic to
a poset of the form
$P_{q_{1}, \ldots, q_{r}}$
described in 
Example~\ref{example_rank_3}.

\item[(ii)] $C(2m) = 2$ for $m \geq 2$ and
the two coatoms in a length $2m$ Sheffer interval
cover exactly the same elements
of rank $2m-2$.

\item[(iii)] $C(2m+1)= h$ is an even positive integer, for $m \geq 2$.
Moreover, the set of $h$ coatoms in a Sheffer interval
of length $2m+1$ partitions into 
$h/2$ pairs,
$\{c_{1},d_{1}\}$, $\{c_{2},d_{2}\}$, $\ldots$,
$\{c_{h/2},d_{h/2}\}$, such that
$c_{i}$ and $d_{i}$ cover the same two elements
of rank $2m-1$.
\end{itemize}
\label{theorem_Sheffer_butterfly}
\end{theorem}
\begin{proof}
Part (i) is immediate since $A(2) \leq C(3)$.
Next we prove (ii).
Let $j = 2m$.
In this case the Euler-Poincar\'e relation
for a Sheffer $j$-interval states:
\begin{equation}
      \sum_{k=1}^{j}
          (-1)^{k}
             \cdot
          \frac{1}{D(k) \cdot 2^{j-k-1}}   
   =
      0        .  
\label{equation_v}
\end{equation}
Use equation~(\ref{equation_v})
in the case of a $(j-2)$-interval to
eliminate the first $j-2$ terms 
in the $j$-interval case of~(\ref{equation_v}),
giving the equality in (ii).
Since $D(j)/(D(j-2) \cdot B(2)) = D(j-1)/(D(j-2) \cdot B(1))$,
the two coatoms in the Sheffer
$j$-interval cover the same elements of rank $j-2$.

Finally, we consider (iii). Assume that $C(j+1) = h$,
where $j = 2m$.
Let $[\hz,y]$ be a Sheffer interval of rank $j+1$.
The number of elements of rank $j$ and of rank $j-1$ 
are both given by $h$.
Moreover each element of rank $j-1$ is covered by
exactly $2$ elements of rank $j$, and by part (ii),
each element of rank $j$ covers
$2$ elements of rank $j-1$.
Hence the order relations between elements of rank
$j-1$ and $j$ are those of 
rank $1$ and $2$ in the poset
$P_{q_{1}, \ldots, q_{r}}$
in
Example~\ref{example_rank_3},
where
$q_{1} + \cdots + q_{r} = h$.

Let $z_{1}, \ldots, z_{q}$ be $q$
coatoms in the Sheffer $(j+1)$-interval
$[\hz,y]$ such that
$z_{i}$ covers $w_{i}$ and $w_{i-1}$,
where we count modulo $q$ in the indices.
That is, $z_{1}$ through $z_{q}$ correspond to the edges in
a $q$-gon and $w_{1}$ through $w_{q}$ to the vertices.
Consider an element $x$ of rank $j-2$
that is covered by $w_{1}$.
The interval $[x,y]$ is isomorphic to $T_{3}$,
that is, the interval has exactly $2$ atoms and
$2$ coatoms.
In this interval the element
$x$ is covered by one more element of rank $j-1$.
Call it $v$.
If the element $v$ does not
correspond
to the elements $w_{2}, \ldots, w_{q}$,
we obtain the contradiction that the interval $[x,y]$
has $4$ coatoms.
If $v$ belongs to 
the elements $w_{2}, \ldots, w_{q}$,
say $w_{i}$,
then the interval $[x,y]$
has the coatoms $z_{1}, z_{2}, z_{i}, z_{i+1}$.
When $q \geq 3$ the set
$\{z_{1}, z_{2}, z_{i}, z_{i+1}\}$
has
at least $3$ members.
Hence the only possibility is that
$q = 2$ and $v = w_{2}$.
Also the coatoms $z_{1}$ and $z_{2}$ 
cover the same elements of rank $j-1$.

We conclude that the only possibility
is that all $q_i$'s are equal to $2$, that is,
$q_{1} = \cdots = q_{r} = 2$. Hence $r = h/2$
and $h$ is an even integer.
Moreover, we also obtain a pairing of the coatoms
such that
the two coatoms in each pair cover
the same elements.
\end{proof}

\begin{theorem}
Let $P$ be an Eulerian Sheffer poset
with the binomial factorial function
satisfying
$B(0) = 1$ and
$B(n) = 2^{n-1}$ for $n \geq 1$
and coatom function $C(n)$.
Then a Sheffer $n$-interval $[\hz,y]$ of
$P$ factors in the rank product as
$[\hz, y] \cong (T_{n-2} \cup \{\hz, \widehat{-1}\}) * Q$,
where
$T_{n-2} \cup \{\hz, \widehat{-1}\}$ denotes the butterfly
interval of rank $n-2$ with two new minimal elements attached
in order,
and $Q$ denotes a poset of rank $n$ such that
\vspace*{-4mm}
\begin{enumerate}
\item[(i)] each element of rank $2$ through $n-1$ in $Q$ is covered
by exactly one element,

\vspace*{-2.5mm}

\item[(ii)] each element of rank $1$ in $Q$ is covered
by exactly two elements,

\vspace*{-2.5mm}

\item[(iii)] each element of even rank $4$ through $2 \lfloor n/2 \rfloor$
in $Q$
covers exactly one element,

\vspace*{-2.5mm}

\item[(iv)] each element of odd rank $k$ from
$5$ through $2 \lfloor n/2 \rfloor + 1$
in $Q$ 
covers exactly $C(k)/2$ elements, and

\vspace*{-2.5mm}

\item[(v)] each $3$-interval $[\hz,x]$ 
in $Q$ is isomorphic to
a poset of the form
$P_{q_{1}, \ldots, q_{r}}$
where $q_{1} + \cdots + q_{r} = C(3)$.
\end{enumerate}
\vspace*{-4mm}
\label{theorem_Sheffer_butterfly_II_second_try}
\end{theorem}
Observe that the poset $Q$ without the minimal element $\hz$
and its atoms forms a tree.
The two posets $Q$ and $T_{n-2} \cup \{\hz, \widehat{-1}\}$
are not Sheffer posets. However, they are
triangular posets.  See the concluding remarks.

\begin{proof_}
{{\bf Proof of Theorem~\ref{theorem_Sheffer_butterfly_II_second_try}:}}
Starting from rank $n-1$
down to rank $3$,
we can partition the elements of
rank $k$ into pairs
using Theorem~\ref{theorem_Sheffer_butterfly}.
To ease notation, partition the remaining ranks
($0$, $1$, $2$ and $n$) into singletons.
This partition
respects the partial order of the interval $[\hz,y]$.
That is, given two blocks $B$ and $C$ such that there exist
two elements $b \in B$ and $c \in C$
so that $b < c$ then for all $b^{\prime} \in B$
and for all $c^{\prime} \in C$ we have that
$b^{\prime} < c^{\prime}$.
Note that this defines a partial order on the blocks.
Denote this poset by $Q$.
It is now straightforward to verify that $Q$ satisfies
the conditions (i) through (v).

To reconstruct the interval $[\hz,y]$ we only have to double
the ranks $3$ through $n-1$. But this is exactly what the rank
product with the poset $T_{n-2} \cup \{\hz, \widehat{-1}\}$ does.
\end{proof_}

\section{Eulerian Sheffer posets with factorial function $B(n) = n!$}
\label{section_Eulerian_Sheffer_with_n_factorial}
\setcounter{equation}{0}

In this section we will classify Eulerian Sheffer posets
that have the factorial function $B(n) = n!$,
that is, every interval $[x,y]$, where $x > \hz$, is
a Boolean algebra.
\begin{theorem}
Let $P$ be an Eulerian Sheffer poset
with binomial factorial function $B(n) = n!$.
Then the
Sheffer factorial function $D(n)$
satisfies one of the following three alternatives:
\begin{itemize}
\item[(i)] $D(n) = 2 \cdot (n-1)!$.
In this case every Sheffer $n$-interval is
of the form $\Sigma^{*}(B_{n-1})$.

\item[(ii)] $D(n) = n!$. In this case the poset
is a binomial poset and hence every Sheffer $n$-interval
is isomorphic to the Boolean algebra $B_{n}$.

\item[(iii)] $D(n) = 2^{n-1} \cdot (n-1)!$.
If we furthermore assume that a Sheffer $n$-interval $[\hz,y]$
is a lattice then the interval $[\hz,y]$ is isomorphic to the cubical 
lattice $C_{n}$.
\end{itemize}
\label{theorem_Sheffer_Boolean}
\end{theorem}
The cubical posets of Example~\ref{example_cubical}
and Example~\ref{example_not_cubical}
demonstrate there is no straightforward
classification of the non-lattice Sheffer intervals
in case (iii) of Theorem~\ref{theorem_Sheffer_Boolean}.
The following examples
further illustrates
Sheffer posets (both finite and infinite)
having the same factorial functions as the cubical
poset.

\begin{example}
{\rm
Let $C_{n}$ be the finite cubical lattice,
that is,
the face lattice of an $(n-1)$-dimensional
cube.
We are going to deform this lattice as follows.
The $1$-skeleton of the cube is a bipartite graph.
Hence the set of atoms $A$ has a natural decomposition
as $A_{1} \cup A_{2}$. Every rank $2$ element 
(edge) covers exactly one atom in each $A_{i}$.
Consider the poset
$$   H_{n} = (C_{n} - A) \cup \left(A_{1} \times \{1,2\} \right)  .  $$
That is, we remove all the atoms and add in
two copies of each atom from $A_{1}$.
Define the cover relations for the new elements
as follows. If $a$ in $A_1$ is covered by $b$
then let $b$ cover both copies $(a,1)$ and $(a,2)$.
The poset $H_{n}$ is a Sheffer poset
with the cubical factorial functions.
}
\label{example_not_cubical_II}
\end{example}

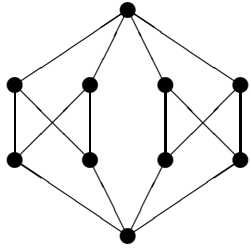
\begin{figure}
\setlength{\unitlength}{1 mm}
\begin{center}
\begin{picture}(30,30)(0,0)
\put(15,0){\circle*{2}}
\put(0,10){\circle*{2}}
\put(10,10){\circle*{2}}
\put(20,10){\circle*{2}}
\put(30,10){\circle*{2}}
\put(0,20){\circle*{2}}
\put(10,20){\circle*{2}}
\put(20,20){\circle*{2}}
\put(30,20){\circle*{2}}
\put(15,30){\circle*{2}}

\put(15,0){\line(3,2){15}}
\put(15,0){\line(1,2){5}}
\put(15,0){\line(-1,2){5}}
\put(15,0){\line(-3,2){15}}

\put(0,10){\line(0,1){10}}
\put(10,10){\line(0,1){10}}
\put(0,10){\line(1,1){10}}
\put(10,10){\line(-1,1){10}}

\put(20,10){\line(0,1){10}}
\put(30,10){\line(0,1){10}}
\put(20,10){\line(1,1){10}}
\put(30,10){\line(-1,1){10}}

\put(0,20){\line(3,2){15}}
\put(10,20){\line(1,2){5}}
\put(20,20){\line(-1,2){5}}
\put(30,20){\line(-3,2){15}}

\end{picture}
\end{center}
\caption{A finite Sheffer poset with the same factorial functions
as the cubical lattice.}
\label{figure_example}
\end{figure}

The poset in Figure~\ref{figure_example}
is the atom deformed cubical lattice $H_{3}$.
This poset is also obtained as length $3$ Sheffer interval
in Example~\ref{example_not_cubical}.

\begin{example}
{\rm
Let $P$ and $Q$ be two Sheffer posets
(finite or infinite)
having the cubical factorial functions
$B(n) = n!$
and
$D(n) = 2^{n-1} \cdot (n-1)!$.
Their diamond product, namely
$P \diamond Q = (P - \{\hz\}) \times (Q - \{\hz\}) \cup \{\hz\}$,
also has the cubical factorial functions.
}
\label{example_not_cubical_III}
\end{example}

\begin{example}
{\rm
As an extension of the previous example,
let $P$ be a Sheffer poset
(finite or infinite)
having the cubical factorial functions.
Then for a set $X$
the poset $(P - \{\hz\})^{X} \cup \{\hz\}$ is
a Sheffer poset with the cubical factorial functions.
The cubical poset (Example~\ref{example_cubical})
is an illustration of this.
}
\label{example_not_cubical_IV}
\end{example}

If we require the extra condition that every finite Sheffer interval
is
a lattice,
we obtain  it is in fact the infinite cubical lattice.
\begin{proposition}
Let $P$ be a finite Sheffer poset of rank $n$ with the cubical
factorial functions
$B(k) = k!$ for $k \leq n-1$
and
$D(k) = 2^{k-1} \cdot (k-1)!$ for $1 \leq k \leq n$.
If $P$ is a lattice then $P$ is
isomorphic to the cubical lattice $C_{n}$.
\label{proposition_cubical_lattice_structure}
\end{proposition}
\begin{proof}
The proof is by induction on the rank $n$ of $P$.
The induction base $n \leq 2$ is straightforward to verify.
Assume true for all posets of rank $n-1$ and consider
a rank $n$ poset $P$.
Using the cubical factorial functions, we know that
the half open interval
$(\hz,\ho]$ contains $3^{n-1}$ elements.
Let $c$ be a coatom in the poset.
The interval $[\hz,c]$ is isomorphic to $C_{n-1}$ by
the induction hypothesis.
Now define a function
$\varphi : (\hz,c] \longrightarrow (\hz,\ho] - (\hz,c]$ 
as follows. For $z$ in $(\hz,c]$
let $\varphi(z)$ be the unique atom in the interval
$[z,\ho]$ that does not belong to the interval $[z,c]$.
The existence and uniqueness follows from the fact the atom
function satisfies $A(k) - A(k-1) = 1$.
Also note that $\varphi(z)$ covers the element $z$.

We next verify the function $\varphi$ is injective. 
If we have $\varphi(z) = \varphi(w)$ then $z$ and $w$ have the same
rank. Also observe that $\varphi(z) \not\leq c$ by the definition
of the function $\varphi$. This contradicts that the interval $[\hz,\ho]$
is a lattice, since $z$ and $w$ have the two upper bounds $\varphi(z)$ and $c$.

The function $\varphi$ also preserves the cover relations.
If $z \coveredby w$
the two-interval $[z,\varphi(w)]$ contains two atoms
which must be $w$ and $\varphi(z)$. Hence
$\varphi(z) \coveredby \varphi(w)$.
Let $\Phi$ be the image of the function $\varphi$.
By a similar argument
the inverse function $\varphi^{-1} : \Phi \longrightarrow (\hz,c]$
also preserves the cover relations.
Thus as posets 
$(\hz,c]$ and $\Phi$ are isomorphic.
Moreover, the disjoint union
$(\hz,c] \cup \Phi$ is an upper order ideal of the 
poset $P$ and has cardinality $2 \cdot 3^{n-2}$.

The poset $P$ has $C(n) = 2n - 2$ coatoms.
One of them is the coatom $c$. Since $c$ covers
$2n-4$ elements there are $2n-4$ coatoms in $\Phi$.
Hence there is a unique coatom
$d$ that does not belong to the upper order ideal $(\hz,c] \cup \Phi$.
Since the interval $[\hz,d]$ is isomorphic to the cubical lattice
$C_{n-1}$ and has $3^{n-2} + 1$ elements, we conclude that
the complement of the upper order ideal is the lower order ideal
$[\hz,d]$.
Thus we have the partition
$(\hz,c] \cup \Phi \cup (\hz,d]$ of $P - \{\hz\}$.

It remains to show that there is a bijective function $\psi : (\hz,d]
\longrightarrow \Phi$ such that $\psi(z)$ covers $z$ and $\psi$ preserves the
cover relation.  Define $\psi : (\hz,d] \longrightarrow (\hz,y] - (\hz,d]
= (\hz,c] \cup \Phi$ 
by letting $\psi(z)$ be the unique atom in the interval
$[z,\ho]$ that does not belong to the interval $[z,d]$.
Observe that if $\psi(z)
\in (\hz,c]$ we obtain that $z < \psi(z) \leq c$, contradicting that
$(\hz,c]$ and $(\hz,d]$ are disjoint. Hence the image of $\psi$ is $\Phi$.
The remaining properties of $\psi$ are proven just like those
for  the function $\varphi$.

Hence $P - \{\hz\}$ is isomorphic to the Cartesian product
of the three element poset
\setlength{\unitlength}{0.4 mm}
\begin{picture}(10,10)(1,0)
\put(0,0){\circle*{2}}
\put(10,0){\circle*{2}}
\put(5,10){\circle*{2}}
\put(0,0){\line(1,2){5}}
\put(10,0){\line(-1,2){5}}
\end{picture}
with $(\hz,c] \cong C_{n-1}$. That is,
the poset is isomorphic to the cubical lattice $C_{n}$.
\end{proof}

\begin{example}
{\rm
Define $\Ccc_{X}^{<\infty}$
to be a subposet of
the cubical poset 
$\Ccc_{X} = P^{X} \cup \{\hz\}$
in Example~\ref{example_cubical},
where
$P$ is the three element poset
\setlength{\unitlength}{0.5 mm}
\begin{picture}(15,10)(-1,0)
\put(0,0){\circle*{2}}
\put(10,0){\circle*{2}}
\put(5,10){\circle*{2}}
\put(0,0){\line(1,2){5}}
\put(10,0){\line(-1,2){5}}
\put(-4,0){${\scriptstyle 0}$}
\put(11,0){${\scriptstyle 1}$}
\put(6,10){${\scriptstyle *}$}
\end{picture},
given by
$$   \Ccc_{X}^{<\infty}
   =
     \{f \in P^{X} \:\: : \:\: |f^{-1}(1)| < \infty\}  \cup \{\hz\} .  $$
That is, for each function $f$ only a finite number of elements of $X$
take on non-zero values.
Since the union of two finite sets is finite it follows
that the join of the two elements is defined. 
It follows that 
$\Ccc_{X}^{<\infty}$ is a lattice.
Observe 
the subposet $\Ccc_{X}^{<\infty}$ 
remains a Sheffer poset with
the cubical factorial functions
$B(n) = n!$ and $D(n) = 2^{n-1} \cdot (n-1)!$.
Call this poset the
{\em infinite cubical lattice}.
}
\label{example_cubical_lattice_indeed}
\end{example}

\begin{theorem}
Let $L$ be an Eulerian Sheffer poset that is also a lattice.
Then $L$ is either isomorphic to $\Bbb_{X}$
where $X$ is the set of atoms of $L$
or $L$ is the infinite cubical lattice
$\Ccc_{X}^{<\infty}$
where $X$ is the set of rank~$2$ elements
of $L$ which are greater than some fixed atom $a$ in $L$.
\end{theorem}
\begin{proof}
Using Theorem~\ref{theorem_Eulerian_binomial_lattice}
we know that the binomial factorial function is
$B(n) = n!$.
Since every Sheffer interval is a lattice
there are only two choices
for the Sheffer factorial function.
The case $D(n) = n!$ is indeed the Boolean algebra
which is the first alternative of the conclusion of the theorem.
Hence let us consider the second choice
$D(n) = 2^{n-1} \cdot (n-1)!$.
Thus every interval $[\hz,y]$ is
a finite cubical lattice.

Let $a$ be an atom of the lattice $L$ and let $X$ be the set of
elements of rank $2$ which cover $a$.  Define the function $\varphi : L
\longrightarrow \Ccc_{X}^{<\infty}$ as follows.
Set $\varphi(\hz) = \hz$.
For $x \in L$ and $x > \hz$ let $y$ be the join of $a$ and~$x$.
Since the interval
$[\hz,y]$ is a finite cubical lattice,
the non-minimal elements of this interval can be encoded
by functions $g : Y \longrightarrow P$, where 
is $P$ is the three element poset 
in 
Example~\ref{example_cubical_lattice_indeed}.
Furthermore we may assume that the set $Y$ is all the elements
in the interval $[a,y]$ that cover $a$.
Without loss of generality,
we may choose the encoding so
that the atom $a$ is the constant function~$0$.

Encode the element $x$ as such a function $g: Y \longrightarrow P$.
Observe that $g$ does not take the value $0$, since that would
contradict that the join of $a$ and $x$ is $y$.
Now define $f : X \longrightarrow P$ by
$$   f(z) = \left\{ \begin{array}{c l}
                 g(z) & \mbox{ if } z \in Y, \\
                  0   & \mbox{ if } z \in X - Y.
                    \end{array} \right. $$
Observe that since $Y$ is a finite set, we know that
$f$ belongs to 
the lattice $\Ccc_{X}^{<\infty}$.
Hence set $\varphi(x)$ to be the function $f$.

The inverse of $\varphi$ is given as follows.
For $f$, a non-zero element of 
the lattice $\Ccc_{X}^{<\infty}$
let the set $Y$ be defined as
$$   Y = \{z \in X \:\: : \:\: f(z) \neq 0\}  .  $$
In the lattice $L$ let the element $y$ be the join
$\bigvee_{z \in Y} z$. Observe that $a \leq y$.
Since the interval $[\hz,y]$ is isomorphic to
the finite cubical lattice $\Ccc_{Y}$,
let $x$ be the unique element corresponding
to the function~$f$ restricted to $Y$.
That is, the inverse of 
$\varphi$ is given by
$\varphi^{-1}(f) = x$.
Moreover let $\varphi^{-1}(\hz) = \hz$.

Observe that both $\varphi$ and $\varphi^{-1}$ are order preserving,
thus proving that 
the lattices $L$ and $\Ccc_{X}^{<\infty}$ are isomorphic.
\end{proof}

Note that it is enough to work with the join operation in this proof,
since a locally finite join semi-lattice with unique minimal element
is a lattice~\cite[Proposition 3.3.1]{Stanley_b}.

We now return to the main issue of classifying
the factorial functions of Eulerian Sheffer posets.
Similar to Lemma~\ref{lemma_A_A}
we have the following lemma.
\begin{lemma}
Let $P$ and $P^{\prime}$ be two Eulerian Sheffer posets
with $B(n) = B^{\prime}(n)$ and
having coatom functions $C(n)$ and $C^{\prime}(n)$ 
which agree
for $n \leq 2m$, where $m \geq 2$.
Then the two following equalities hold:
\begin{equation}
        \frac{1}{C(2m+1)}
     \cdot
        \left(
          1 - \frac{2}{C(2m+2)}
        \right) 
    =
        \frac{1}{C^{\prime}(2m+1)}
     \cdot
        \left(
          1 - \frac{2}{C^{\prime}(2m+2)}
        \right), 
\label{equation_C_C_1}
\end{equation}
and
\begin{eqnarray}
&&
        \frac{1}{C(2m+1)}
     \cdot
        \left(
          \frac{1}{B(3)} - \frac{1}{C(2m+2)}
     \cdot
        \left(
          \frac{1}{2} - \frac{1}{C(2m+3)}
     \cdot
        \left(
          1 - \frac{2}{C(2m+4)}
        \right) 
        \right) 
        \right)      \nonumber \\
& = &
        \frac{1}{C^{\prime}(2m+1)}
     \cdot
        \left(
          \frac{1}{B(3)} - \frac{1}{C^{\prime}(2m+2)}
     \cdot
        \left(
          \frac{1}{2} - \frac{1}{C^{\prime}(2m+3)}
     \cdot
        \left(
          1 - \frac{2}{C^{\prime}(2m+4)}
        \right) 
        \right) 
        \right)     . 
\label{equation_C_C_2} 
\end{eqnarray}
\label{lemma_C_C}
\end{lemma}

Similar to 
Corollary~\ref{corollary_A_A_bound}
we have the following result.
\begin{corollary}
Let $P$ and $P^{\prime}$ be two Eulerian Sheffer posets 
satisfying the same conditions as in Lemma~\ref{lemma_C_C}.
Assume furthermore that there is a lower and an upper bound
for $C^{\prime}(2m+2)$ of the form
$L \leq C^{\prime}(2m+2) < U$. 
Let $x$ be the left-hand side of equation~(\ref{equation_C_C_1}).
Then we obtain 
a lower and an upper bound for
$C^{\prime}(2m+1)$, namely
\begin{equation}
     \frac{1}{x} \cdot \left(1 - \frac{2}{L}\right)
   \leq
     C^{\prime}(2m+1)
   <
     \frac{1}{x} \cdot \left(1 - \frac{2}{U}\right) .  
\label{equation_C_C_1_bound}
\end{equation}
Similarly,
let $z$ be the left-hand side of equation~(\ref{equation_C_C_2})
and let
$$  y = 
         \frac{1}{2}
       -
           C^{\prime}(2m+2)
         \cdot
           \left(
             \frac{1}{B(3)} - C^{\prime}(2m+1) \cdot z
           \right)
       .
$$
Then the lower and upper bound
$L \leq C^{\prime}(2m+4) < U$
implies
\begin{equation}
     \frac{1}{y} \cdot \left(1 - \frac{2}{L}\right)
   \leq
     C^{\prime}(2m+3)
   <
     \frac{1}{y} \cdot \left(1 - \frac{2}{U}\right) .  
\label{equation_C_C_2_bound}
\end{equation}
\label{corollary_C_C_bound}
\end{corollary}
Both bounds can be improved by using that
$C^{\prime}(2m+1)$ and
$C^{\prime}(2m+3)$
are integers.

The proof of 
the main result of this section,
Theorem~\ref{theorem_Sheffer_Boolean},
is broken down into
four propositions,
namely
Propositions~\ref{proposition_i},
\ref{proposition_ii},
\ref{proposition_iii}
and~\ref{proposition_iv}.
The proof of each proposition
branches into several cases
and one has to show that these cases cannot occur.
The main tool to exclude these possibilities
are Fact~\ref{fact_a_b}
and the bounds in Corollary~\ref{corollary_C_C_bound}.
In one case we use Proposition~\ref{proposition_special_case}.

\begin{proposition}
Let $P^{\prime}$ be an Eulerian Sheffer poset
with factorial functions
satisfying
$B^{\prime}(n) = n!$
and
$D^{\prime}(3) = 4$.
Then the Sheffer factorial function
is given by $D^{\prime}(n) = 2 \cdot (n-1)!$.
\label{proposition_i}
\end{proposition}
\begin{proof}
Let $P$ be the poset $\Sigma^{*}(\Bbb)$
with the coatom function
$C(n) = n-1$ for $n \geq 3$.

Assume that the coatom functions $C$ and $C^{\prime}$
agree for $n \leq 2m = j$.
Then the left-hand side of equation~(\ref{equation_C_C_1})
is given by $(j-1)/(j (j+1))$.
The bounds on $C^{\prime}(j+2)$
are $j+1 = A(j+1) \leq C^{\prime}(j+2) < \infty$.
Now from~(\ref{equation_C_C_1_bound})
we have
$$ j \leq C^{\prime}(j+1) < j+2+\frac{2}{j-1}  .  $$
Since $j \geq 4$ we have three cases
$C^{\prime}(j+1) = j, j+1, j+2$.
\begin{itemize}
\item[(a)] The case $C^{\prime}(j+1) = j+1$.
Consider a rank $j+1$ Sheffer interval.
It has $D^{\prime}(j+1)/B(j)$ atoms.
However
$D^{\prime}(j+1)/B(j) = C^{\prime}(j+1) \cdot D(j)/B(j) = 
(j+1) \cdot 2 \cdot (j-1)!/j! = 
2 \cdot (j+1)/j = 
(2 \cdot m+1)/m = 2 + 1/m$,
which is not an integer for $m \geq 2$.

\item[(b.i)] The case $C^{\prime}(j+1) = j+2$ and we assume
$j \geq 6$.
This is done similarly as the previous
case.
The number of atoms is given by
$D^{\prime}(j+1)/B(j) = 2 + 2/m$,
which is not an integer for $m \geq 3$.

\item[(b.ii)] The case $C^{\prime}(j+1) = j+2$
when $j=4$, that is,
$C^{\prime}(5) = 6$ and
$C^{\prime}(6) = 20$.
Equation~(\ref{equation_C_C_2})
implies 
$1/C^{\prime}(7) 
  \cdot
 (1 - 2/C^{\prime}(8))
    =
      - 5/42$,
which does not have any positive integer solutions.
\end{itemize}
The remaining case is $C^{\prime}(j+1) = j$
which implies $C^{\prime}(j+2) = j+1$.
Hence the two coatom functions~$C$ and $C^{\prime}$ are equal.
\end{proof}

\begin{lemma}
Let $P$ be a rank $n$ finite Eulerian Sheffer poset
with factorial functions
$B(k) = k!$ for $k \leq n-1$
and
$D(k) = 2 \cdot (k-1)!$ for $2 \leq k \leq n$.
Then the poset $P$ is isomorphic to $\Sigma^{*}(B_{n-1})$.
\label{lemma_i_structure}
\end{lemma}
\begin{proof}
Observe that $P$ has $D(n)/B(n-1) = 2$ atoms.
Denote them by $a_{1}$ and $a_{2}$.
Also note that every element of rank $2$ in $P$
covers both atoms. Finally, since the interval
$[a_{i},\ho]$ is isomorphic to~$B_{n-1}$, we obtain
that $P$ is isomorphic to $\Sigma^{*}(B_{n-1})$.
\end{proof}

\begin{proposition}
Let $P^{\prime}$ be an Eulerian Sheffer poset
with factorial functions
satisfying
$B^{\prime}(n) = n!$
and
$D^{\prime}(3) = 6$.
Then the factorial function
is given by $D^{\prime}(n) = n!$.
\label{proposition_ii}
\end{proposition}
\begin{proof}
Let $P$ be the infinite Boolean algebra $\Bbb$
with coatom function $C(n) = n$.

Assume that $C(n)$ and $C^{\prime}(n)$ are equal
for all $n \leq 2m = j$. Now we have the bound
$j+1 = A(j+1) \leq C^{\prime}(j+2) < \infty$.
Corollary~\ref{corollary_C_C_bound}
implies
$j+1-2/j \leq C^{\prime}(j+1) < j+3+2/j$.
That is, we have
$j+1 \leq C^{\prime}(j+1) \leq j+3$.

\begin{itemize}
\item[(a)]
$C^{\prime}(j+1) = j+2$.
This case is ruled out by
Proposition~\ref{proposition_special_case}
since a finite Sheffer poset of rank~$j+1$
having these factorial functions would be a finite
binomial poset.

\item[(b)]
$C^{\prime}(j+1) = j+3$.
Equation~(\ref{equation_C_C_1})
implies
$C^{\prime}(j+2) = (j+1) \cdot (j+2)$.
Now equation~(\ref{equation_C_C_2})
states
$1/C^{\prime}(j+3) \cdot (1-2/C^{\prime}(j+4))
  =  - (j^2-4)/(6 \cdot (j+4))$, which is negative for $j \geq 4$.
\end{itemize}
The only remaining case is
$C^{\prime}(j+1) = j+1$ which implies
$C^{\prime}(j+2) = j+2$. Hence
the two coatom functions $C$ and $C^{\prime}$ are identical.
\end{proof}

\begin{proposition}
Let $P^{\prime}$ be an Eulerian Sheffer poset
with factorial functions
satisfying
$B^{\prime}(n) = n!$
and
$D^{\prime}(3) = 8$.
Then the factorial function 
is given by $D^{\prime}(n) = 2^{n-1} \cdot (n-1)!$.
\label{proposition_iii}
\end{proposition}
\begin{proof}
Let $P$ be the cubical lattice with
coatom function $C(n) = 2 \cdot (n-1)$
and factorial function $D(n) = 2^{n-1} \cdot (n-1)!$.
Assume that 
the coatom functions
$C$ and $C^{\prime}$ agree up to
$2m = j$. 
Using 
Corollary~\ref{corollary_C_C_bound}
with the bounds
$j+1 = A(j+1) \leq C^{\prime}(j+2) < \infty$
we obtain
$2j-2 \leq C^{\prime}(j+1) \leq 2j+1$.

The two bounds
$j+2 \leq C^{\prime}(j+3) < \infty$
and
$j+3 \leq C^{\prime}(j+4) < \infty$
give the bound
\begin{equation}
     0
  <
     \frac{1}{C^{\prime}(j+3)}
        \cdot
     \left(1 - \frac{2}{C^{\prime}(j+4)}\right)
  <
     \frac{1}{j+2}   .  
\label{equation_iii_four}
\end{equation}

Consider now the cases:
\begin{itemize}
\item[(a)] $C^{\prime}(j+1) = 2j-2$.
Now equation~(\ref{equation_C_C_1}) implies
$C^{\prime}(j+2) = j+1$.
Equation~(\ref{equation_C_C_2}) states
$1/C^{\prime}(j+3)  \cdot (1 - 2/C^{\prime}(j+4))
 = (j+7)/(12 \cdot (j+3))$.

\item[(a.i)] When $j \geq 8$ we have that
$1/C^{\prime}(j+3)  \cdot (1 - 2/C^{\prime}(j+4))
 = (j+7)/(12 \cdot (j+3)) > 1/(j+2)$,
contradicting inequality~(\ref{equation_iii_four}).

\item[(a.ii)] $j = 4$.
Then we have
$C^{\prime}(5) = 6$ and
$C^{\prime}(6) = 5$.
Now we have the identity
$1/C^{\prime}(7)  \cdot (1 - 2/C^{\prime}(8))
 = 11/84$.
Hence the inequality
$7 \leq C^{\prime}(8) < \infty$ implies
$60/11 \leq C^{\prime}(7) < 84/11$.
That is, 
$6 \leq C^{\prime}(7) \leq 7$.
However, $C^{\prime}(7) = 6$ implies
$C^{\prime}(8) = 28/3$, not an integer.
Hence the only possible case is
$C^{\prime}(7) = 7$.

The number of elements of rank $5$ in a rank $7$
Sheffer interval is given by
$D^{\prime}(7)/(D^{\prime}(5) \cdot B(2)) =
 C^{\prime}(7) \cdot C^{\prime}(6)/2 =
 7 \cdot 5 /2$, which is not an integer.

\item[(a.iii)] $j = 6$.
Then we have
$C^{\prime}(7) = 10$ and
$C^{\prime}(8) = 7$.
The numbers of atoms in a Sheffer interval of rank $7$
is given by
$D^{\prime}(7)/B(6) = 
 C^{\prime}(7) \cdot D^{\prime}(6)/B(6) = 
 10 \cdot 2^{5} \cdot 5! / 6! = 
 5 \cdot 2^{5}/3$
which is not an integer.

\item[(b)] The case when
$C^{\prime}(j+1) = 2j-1$.
Now equation~(\ref{equation_C_C_1}) implies
$C^{\prime}(j+2) = (4j+4)/3$.
Equation~(\ref{equation_C_C_2}) implies
$1/C^{\prime}(j+3)  \cdot (1 - 2/C^{\prime}(j+4))
 =  (j+10)/(18 \cdot (j+3))$.
Also since $(4j+4)/3$ is an integer, we have
the congruence condition $j \equiv 2 \bmod 6$.

\item[(b.i)] $j \geq 14$.
Now $1/C^{\prime}(j+3)  \cdot (1 - 2/C^{\prime}(j+4))
 =  (j+10)/(18 \cdot (j+3)) > 1/(j+2)$ as $j \geq 14$.

\item[(b.ii)] $j = 8$.
Then we have
$C^{\prime}(9) = 15$.
Now equation~(\ref{equation_C_C_1}) implies
$C^{\prime}(10) = 12$.
Equation~(\ref{equation_C_C_2}) states
$1/C^{\prime}(11)  \cdot (1 - 2/C^{\prime}(12)) =  1/11$.
The bounds
$11 \leq C^{\prime}(12) < \infty$
imply
$9 \leq C^{\prime}(11) < 11$

\item[(b.ii.1)]
$C^{\prime}(11)  = 9$ which implies $C^{\prime}(12) = 11$.
The number of elements of rank $10$ in a 
Sheffer
interval of rank $12$ is given by
$C^{\prime}(12) \cdot C^{\prime}(11)/2 = 99/2$.
Hence this case is excluded.

\item[(b.ii.2)]
$C^{\prime}(11)  = 10$ which implies that $C^{\prime}(12) = 22$.
Now the Euler-Poincar\'e relation on rank $14$ Sheffer
interval implies that
$C^{\prime}(13) = -39/4 \cdot (1 - 2/C^{\prime}(14))$
which has no positive integer solutions.

\item[(c)] The case $C^{\prime}(j+1) = 2j+1$.
Equation~(\ref{equation_C_C_1}) implies
$C^{\prime}(j+2) = 4j+4$.
Equation~(\ref{equation_C_C_2}) implies
$1/C^{\prime}(j+3)  \cdot (1 - 2/C^{\prime}(j+4))
 = - (j-2)/(6 \cdot (j+3))$ which
is negative for $j \geq 4$.
\end{itemize}
The only remaining case is
$C^{\prime}(j+1) = 2j$ which implies
$C^{\prime}(j+2) = 2j+2$. Thus we conclude that
the coatom functions $C$ and $C^{\prime}$ are equal.
\end{proof}

\begin{proposition}
There is no Eulerian Sheffer poset
with factorial functions
$B^{\prime}(n) = n!$
and
$D^{\prime}(3) = 10$.
\label{proposition_iv}
\end{proposition}
\begin{proof}
The Euler-Poincar\'e relation implies that
$C^{\prime}(4) = 12$.
The Euler-Poincar\'e relation on
a Sheffer $6$-interval implies that $C^{\prime}(6) = 2$,
which contradicts $C^{\prime}(6) \geq A^{\prime}(5)$.
\end{proof}

\begin{proof_}{{\bf Proof of Theorem~\ref{theorem_Sheffer_Boolean}:}}
The Euler-Poincar\'e relation for
a Sheffer $4$-interval states
$$         1
        -
           \frac{2}{C(4)}
   =
      \frac{C(3)}{6}    .   $$
Hence $C(3) < 6$, giving
the four possibilities
$C(3) = 2,3,4,5$.
They are addressed in
the four
Propositions~\ref{proposition_i},
\ref{proposition_ii},
\ref{proposition_iii}
and~\ref{proposition_iv}.
Similarly, the structure results
are proved in
Lemma~\ref{lemma_i_structure} and
Propositions~\ref{proposition_Boolean_structure}
and~\ref{proposition_cubical_lattice_structure}.
\end{proof_}

\section{Concluding remarks}
\setcounter{equation}{0}

An interesting research project is to classify the factorial functions
of finite Eulerian binomial posets and finite Eulerian Sheffer posets.
Two examples of finite Sheffer posets are the
face lattices of the dodecahedron and the four-dimensional
regular polytope known as the $120$-cell.
In Propositions~\ref{proposition_Boolean},
\ref{proposition_butterfly},
\ref{proposition_i},
\ref{proposition_ii},
\ref{proposition_iii}
and~\ref{proposition_iv}
many finite possibilities for the factorial functions
were excluded since there was no
possibility to extend the factorial function to
higher ranks.
A first step in this classification is to consider these cases.

Also note the following lemma,
the proof of which follows directly from
Proposition~\ref{proposition_only_even_ranks}.
\begin{lemma}
Let $P$ be an Eulerian finite binomial (Sheffer) poset of odd rank $n$.
Let $Q$ be the poset obtained by taking $k$ disjoint copies of $P$
and identifying the minimal, respectively, maximal elements.
Then Q is an Eulerian finite binomial (Sheffer) poset.
The only value of the factorial function(s) that changes is
the one that enumerates the maximal chains, namely,
$B_{Q}(n) = k \cdot B_{P}(n)$ in the binomial case,
and $D_{Q}(n) = k \cdot D_{P}(n)$ in the Sheffer case.
\end{lemma}

A larger class of posets to consider are
the triangular posets~\cite{D_R_S}.
A poset is {\em triangular} if every interval
$[x,y]$, where $x$ has rank $n$ and $y$ has rank $m$, has
$B(n,m)$ maximal chains.
Both binomial and Sheffer posets are triangular.
A non-trivial Eulerian example
of a finite triangular poset is the
face lattice of the $4$-dimensional
regular polytope known as the $24$-cell.
Can the factorial function $B(n,m)$
be classified for Eulerian triangular posets?

Classifying 
finite Eulerian Sheffer posets
seems to be hard
as seen from the multitude of examples having
the cubical factorial functions.
We leave the reader with three examples
of Sheffer posets having the same factorial
functions as the face lattice of the dodecahedron,
each of which is not isomorphic to 
this face lattice.

\begin{example}
{\rm
{\em An Eulerian finite Sheffer poset with the same
factorial functions as the face lattice of the dodecahedron.}
For an $n$-gon define a $CW$-complex $X_{n}$ as follows.
First take the antiprism of the $n$-gon.
We then have a $CW$-complex consisting of two $n$-gons
and $2n$ triangles. Note that at every vertex three triangles
and one $n$-gon meet. Now subdivide each of the two $n$-gons
by placing a vertex in each $n$-gon and attaching this vertex by
$n$ new edges to the $n$ vertices of the $n$-gon.
Let this be the $CW$-complex $X_{n}$.

Observe that $X_{n}$ consists of $2n+2$ vertices,
$6n$ edges and $4n$ triangles. Moreover, at $2n$ of the vertices
$5$ triangles meet. At the other two vertices $n$ triangles
meet. Label these two vertices $a$ and~$b$.
Also note that $X_{5}$ is the boundary complex of an icosahedron.
Observe for $n \geq 3$ that $X_{n}$ is a simplicial complex.
However, for $n=2$ it is necessary to view $X_{2}$ as a $CW$-complex.


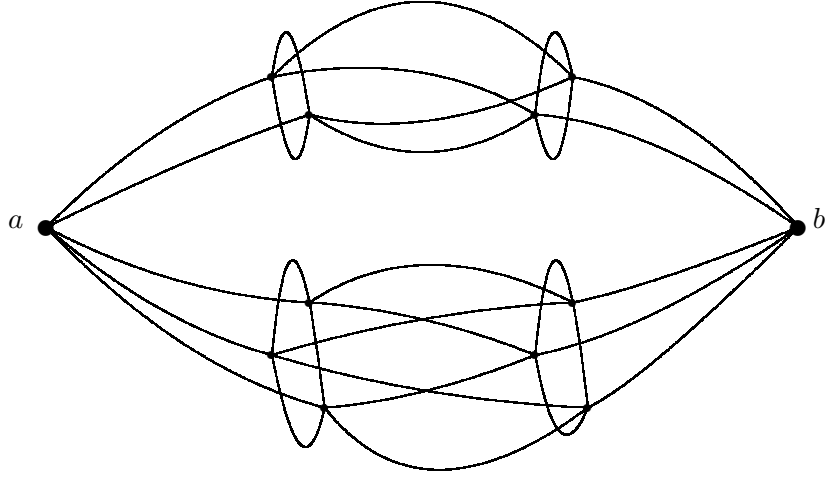
\begin{figure}
\setlength{\unitlength}{1mm}
\begin{center}
\begin{picture}(100,63)(0,17)

\put(0,50){\circle*{2}}
\put(-5,50){$a$}

\put(30,70){\circle*{1}}
\put(35,65){\circle*{1}}

\put(70,70){\circle*{1}}
\put(65,65){\circle*{1}}

\put(100,50){\circle*{2}}
\put(102,50){$b$}

\thinlines


\qbezier(0,50)(15,65)(30,70)
\qbezier(0,50)(20,60)(35,65)


\qbezier(30,70)(32,84)(35,65)
\qbezier(30,70)(33,51)(35,65)

\qbezier(30,70)(50,74)(65,65)
\qbezier(30,70)(50,90)(70,70)
\qbezier(35,65)(50,55)(65,65)
\qbezier(35,65)(50,61)(70,70)

\qbezier(65,65)(67,84)(70,70)
\qbezier(65,65)(68,51)(70,70)

\qbezier(65,65)(80,64)(100,50)
\qbezier(70,70)(85,67)(100,50)

\put(30,33){\circle*{1}}
\put(35,40){\circle*{1}}
\put(37,26){\circle*{1}}

\put(65,33){\circle*{1}}
\put(70,40){\circle*{1}}
\put(72,26){\circle*{1}}


\qbezier(0,50)(15,37)(30,33)
\qbezier(0,50)(17,41)(35,40)
\qbezier(0,50)(18,31)(37,26)

\qbezier(30,33)(32,54)(35,40)
\qbezier(30,33)(34,12.75)(37,26)
\qbezier(35,40)(36,34)(37,26)

\qbezier(35,40)(50,39)(65,33)
\qbezier(37,26)(50,27)(65,33)
\qbezier(30,33)(50,39)(70,40)
\qbezier(30,33)(50,27)(72,26)
\qbezier(35,40)(50,50)(70,40)
\qbezier(37,26)(50,9.5)(72,26)

\qbezier(65,33)(67,54)(70,40)
\qbezier(65,33)(68,16.25)(72,26)
\qbezier(72,26)(71,35)(70,40)

\qbezier(65,33)(80,36)(100,50)
\qbezier(70,40)(80,42)(100,50)
\qbezier(72,26)(80,30)(100,50)



\end{picture}
\end{center}
\caption{The $CW$-complex obtained by
joining the complexes $X_{2}$ and $X_{3}$ at the vertices
$a$ and $b$.}
\label{figure_X_2_X_3}
\end{figure}


Construct a $CW$-complex $Y$ by taking $X_{2}$ and $X_{3}$ and
identifying the vertices labeled $a$ and
identifying the vertices labeled $b$.
See Figure~\ref{figure_X_2_X_3}.
The dual of the face poset of $Y$ is an Eulerian Sheffer poset
with factorial functions agreeing with the face lattice of a dodecahedron.
}
\label{example_dodecahedron_one}
\end{example}

\begin{example}
{\rm
For $1 \leq i \leq 3$
let $Z_{i}$ be the boundary
of a $3$-dimensional simplex
with vertices $z_{i,1}$, $z_{i,2}$, $z_{i,3}$ and $z_{i,4}$.
Similarly,
for $1 \leq j \leq 4$
let $W_{j}$ be the spherical $CW$-complex consisting
of two triangles sharing the three edges.
Call the vertices $w_{1,j}$, $w_{2,j}$ and $w_{3,j}$.
Now identify vertex $z_{i,j}$ with~$w_{i,j}$.
We then have a $CW$-complex 
that has $12$ vertices,
$3 \cdot 6 + 4 \cdot 3 = 30$ edges
and
$3 \cdot 4 + 4 \cdot 2 = 20$ triangles.
Observe that the vertex figure of every vertex
is the disjoint union of a $2$-gon and a triangle.
Thus the dual of the face poset is Sheffer poset
with the same factorial functions as the face lattice of a dodecahedron.
In fact, one may obtain several of these
$CW$-complexes by choosing different identifications
between the two classes of vertices.
}
\label{example_dodecahedron_two}
\end{example}

\begin{example}
{\rm
A third example is formed by taking two $X_{2}$'s from
Example~\ref{example_dodecahedron_one}
and
the boundary of one $3$-dimensional simplex, $Z$, from
Example~\ref{example_dodecahedron_two}
and identifying vertices $a_1$, $a_2$, $b_1$ and $b_2$
with the vertices of the simplex.
}
\end{example}

A different proof of
Proposition~\ref{proposition_Boolean_structure} may
be given using the following
result of Stanley.  A graded finite poset $P$ is a Boolean algebra if
every $3$-interval is a Boolean algebra and for every interval $[x,y]$
of rank of least $4$ the open interval $(x,y)$ is connected.
See~\cite[Lemma~8]{Grabiner}.  Hence it is natural to ask 
if one can extend this result to cubical lattices.
That is, a graded
finite poset $P$ is a cubical lattice if every $3$-interval
$[x,y]$, where $x > \hz$, is a Boolean algebra, every $3$-interval
$[\hz,y]$ is the face lattice of a square, and for every interval
$[x,y]$ of rank of least $4$ the open interval $(x,y)$ is connected.

One may drop the Eulerian condition and ask
to characterize Sheffer posets which are lattices.
The lattice-theoretic
techniques of Farley and Schmidt
may be useful~\cite{Farley_Schmidt}.

Finally, there are long-standing open questions regarding binomial
posets.  One such question asked whether there exist two binomial posets
having the same factorial function but non-isomorphic intervals.  This
question was very recently settled by J\"orgen Backelin~\cite{Backelin}.
However, it is still unknown if there is a binomial
poset having the atom function $A(n) = F_{n}$, the $n$th Fibonacci
number.  See Exercise~78b, Chapter~3 in~\cite{Stanley_b}.

\section*{Acknowledgements}

The first author was partially supported by
National Science Foundation grant 0200624
and by 
a
University of Kentucky College of Arts~\& Sciences 
Faculty Research Fellowship.
The second author was partially supported by
a 
University of Kentucky College of Arts~\& Sciences Research Grant.
Both authors thank G\'abor Hetyei
for inspiring them to
study Eulerian binomial posets,
the Banff International Research Station where some of the
ideas for this paper were
developed,
and the Mittag-Leffler Institute where this paper was completed.
Both authors gratefully acknowledge the careful and thoughtful comments made
by one of the anonymous referees.


\newcommand{\journal}[6]{{\sc #1,} #2, {\it #3} {\bf #4} (#5), #6.}
\newcommand{\book}[4]{{\sc #1,} ``#2,'' #3, #4.}
\newcommand{\thesis}[4]{{\sc #1,} ``#2,'' Thesis, #3, #4.}
\newcommand{\springer}[4]{{\sc #1,} ``#2,'' Lecture Notes in Math.,
                          Vol.\ #3, Springer-Verlag, Berlin, #4.}
\newcommand{\preprint}[3]{{\sc #1,} #2, preprint #3.}
\newcommand{\appear}[3]{{\sc #1,} #2, to appear in {\it #3}.}
\newcommand{\JCTA}{J.\ Combin.\ Theory Ser.\ A}
\newcommand{\arXiv}[4]{{\sc #1,} #2, arXiv: {\tt #3}, #4.}


{\em
\noindent
R.\ Ehrenborg,
Department of Mathematics,
University of Kentucky,
Lexington, KY 40506, \newline
{\tt jrge@ms.uky.edu}
}

{\em
\noindent
M.\ Readdy,
Department of Mathematics,
University of Kentucky,
Lexington, KY 40506, \newline
{\tt readdy@ms.uky.edu}
}


\begin{thebibliography}{99}



\bibitem{Backelin}
\arXiv{J.\ Backelin}
        {Binomial posets with non-isomorphic intervals}
        {math.CO/0508397}
        {22 August 2005}



\bibitem{Bayer_Hetyei}
\journal{M.\ M.\ Bayer and G.\ Hetyei}
        {Flag vectors of Eulerian partially ordered sets}
        {European J.\ Combin.}
        {22}{2001}{5--26}




\bibitem{Billera_Liu}
\journal{L.\ J.\ Billera and N.\ Liu}
        {Noncommutative enumeration in graded posets}
        {J.\ Algebraic Combin.}
        {12}{2000}{7--24}

\bibitem{Bjorner_Welker}
\appear{A.\ Bjorner, V.\ Welker}
       {Segre and Rees products of posets,
        with ring-theoretic applications}
       {J.\ Pure Appl.\ Algebra}


\bibitem{D_R_S}
{\sc P.\ Doubilet, G.-C.\ Rota and R.\ Stanley,}
{\it On the foundation of combinatorial theory (VI).
     The idea of generating functions,}
in Sixth Berkeley Symp.\ on Math.\ Stat.\ and Prob.,
vol. 2: Probability Theory, Univ. of California
(1972), pp. 267--318.



\bibitem{Ehrenborg_Readdy}
\journal{R.\ Ehrenborg and M.\ Readdy}
        {Sheffer posets and $r$-signed permutations}
        {Annales des Sciences Math\'ematiques du Qu\'ebec}
        {19}{1995}{173--196}

\bibitem{Ehrenborg_Readdy_homology}
\journal{R.\ Ehrenborg and M.\ Readdy}
        {Homology of Newtonian coalgebras}
        {European J.\ Combin.}
        {23}{2002}{919--927}


\bibitem{Farley_Schmidt}
\journal{J.\ Farley and S.\ Schmidt}
        {Posets that locally resemble distributive lattices}
        {\JCTA}
        {92}{2000}{119--137}


\bibitem{Grabiner}
\journal{D.\ J.\ Grabiner}
        {Posets in which every interval is a product of chains,
         and natural local actions of the symmetric group}
        {Discrete Math.}
        {199}{1999}{77--84}


\bibitem{Reiner}
\journal{V.\ Reiner}
        {Upper binomial posets and signed permutation statistics}
        {European J.\ Combin.}
        {14}{1993}{581--588}


 
\bibitem{Stanley_a}
\journal{R.\ Stanley}
        {Binomial posets, M\"obius inversion, and permutation enumeration}
        {\JCTA}
        {20}{1976}{336--356}
 

\bibitem{Stanley_b}
\book{R.\ Stanley}
     {Enumerative Combinatorics, Vol. I}
     {Wadsworth and Brooks/Cole, Pacific Grove}
     {1986}




\end{thebibliography}
\end{document}